\documentclass[journal]{IEEEtran}
\makeatletter
\def\endthebibliography{%
	\def\@noitemerr{\@latex@warning{Empty `thebibliography' environment}}%
	\endlist
}
\makeatother
\usepackage{array}
\usepackage{amsmath,amsfonts,amssymb,amsthm}
\usepackage{cite}
\usepackage{graphicx}
\usepackage{mathtools}
\usepackage{epstopdf}
\usepackage{enumitem}
\usepackage[utf8]{inputenc}
\usepackage[english]{babel} 
\usepackage{algorithm, algorithmic}
\usepackage{stfloats}

\usepackage{caption}
\usepackage{subcaption}

\newtheorem*{remark}{Remark}  
\DeclareGraphicsExtensions{.eps}
\newcommand\norm[1]{\left\lVert#1\right\rVert}

\ifCLASSINFOpdf
\else
\fi
\begin{document}
%
\title{A Stochastic Quasi-Newton Method for Large-Scale Nonconvex Optimization with Applications}
%
%
%

\author{H.~Chen,
	~H. C.~Wu,~\IEEEmembership{Member,~IEEE},	
	~S. C. Chan,~\IEEEmembership{Member,~IEEE},	
	~W. H.~Lam,~\IEEEmembership{Senior Member,~IEEE}
	\thanks{The authors are with the Department of Electrical and Electronic Engineering, The University of Hong Kong, Hong Kong, China (e-mail: hmchen@eee.hku.hk;
		andrewhcwu@eee.hku.hk;
		scchan@eee.hku.hk;  
		whlam@eee.hku.hk;).}
}

%
%

\markboth{ IEEE TRANSACTIONS ON NEURAL NETWORKS AND LEARNING SYSTEMS, Dec.~2019}%
{Shell \MakeLowercase{\textit{et al.}}: Bare Demo of IEEEtran.cls for IEEE Journals}
%



\maketitle

\begin{abstract}
Ensuring the positive definiteness and avoiding ill-conditioning of the Hessian update in the stochastic Broyden–Fletcher–Goldfarb–Shanno (BFGS) method are significant in solving nonconvex problems. This paper proposes a novel stochastic version of damped and regularized BFGS method for addressing the above problems. While the proposed regularized strategy helps to prevent the BFGS matrix from being close to singularity, the new damped parameter further ensures positivity of the product of correction pairs. To alleviate the computational cost of the stochastic LBFGS updates, and to improve its robustness, the curvature information is updated using the averaged iterate at spaced intervals. The effectiveness of the proposed method is evaluated through the logistic regression and Bayesian logistic regression problems in machine learning. Numerical experiments are conducted by using both synthetic dataset and several real datasets. The results show that the proposed method generally outperforms the stochastic damped limited memory BFGS (SdLBFGS) method. In particular, for problems with small sample sizes, our method has shown superior performance and is capable of mitigating ill-conditioned problems. Furthermore, our method is more robust to the variations of the batch size and memory size than the SdLBFGS method.
\end{abstract}


%
\IEEEpeerreviewmaketitle

\section{Introduction}
%
%
%
%
\IEEEPARstart{S}{tochastic}  optimization algorithms have been extensively studied over decades and can be traced back to the epochal work \cite{RobbinsHerbert1951ASAM}, which have been widely employed in different areas, e.g., machine learning \cite{ZoubinGhahramani2015Pmla,BleiDavidM.2017VIAR,HoffmanSVI,JRAFATI,SCA}, power systems \cite{SHUANG}, wireless communication \cite{DSOCS,ARibeiro,PSI}, and bioinformatics \cite{MCHEN}. In particular, the classical stochastic approximation (SA) of the exact gradient, also known as stochastic gradient descent (SGD), has been widely applied to these stochastic optimization problems, where the gradient information is employed in finding the search direction. However, in many applications, the exact gradient depends on certain random variables with unknown distributions and thus is difficult to evaluate explicitly. Furthermore, in many applications with extremely massive data samples, the exact gradient of the objective function is rather expensive to compute. In SGD, an unbiased estimator of the gradient is derived using a mini-batch of data points randomly sampled from the full dataset. This substantially reduces the computational cost.

In the theoretical aspect,  SGD algorithm has been widely used in the problems with the assumption that the objective function $f(\cdot)$  is twice continuously differentiable and strongly convex. In particular, \cite{RobustSA} has proposed a robust mirror descent SA algorithm, which is also applicable to general convex objective functions.  Recently, there has been an increasing interest in SA based algorithms for solving nonconvex stochastic optimization problems  \cite{SQN,StoFirZeroSP,SBMD}. Specifically,  \cite{SBMD} has investigated a stochastic block mirror descent method to solve large scale nonconvex optimization problems with high dimensional optimization variables.  \cite{StoFirZeroSP} has studied a framework of randomized stochastic gradient (RSG) methods by randomly selecting a solution from the previous iterates. The Monte Carlo integration has been adopted for the stochastic search direction  \cite{VBISS,YiSun2009Ssut}. Moreover, the control variate technique \cite{VBISS} is proposed to reduce the variance of the SA.

In the deterministic optimization settings, quasi-Newton or Newton methods can achieve higher accuracy and faster convergence by utilizing the second-order information  \cite{SQN,OPREVIEW}. For the stochastic regime, stochastic quasi-Newton’s methods (SQN) have been extensively studied in  \cite{SBBFGS,SGD-QN,MokhtariAryan2014RRSB,SQN,IQN,RHBYRD,RHBYRD2,OPREVIEW,Powell1978,oLBFGS2,AJA}. In particular,  \cite{oLBFGS2} has developed a stochastic variable-metric method with subsampled gradients. In \cite{SGD-QN}, a SGD-QN scheme has been proposed in which the diagonal elements of the Hessian matrix are approximated to rescale the SGD. Since it only involves scalar computation, the method is quite efficient. It should be noted that direct application of the deterministic quasi-Newton methods brings noisy curvature approximation and thus affects the robustness of the iteration  \cite{RHBYRD}. In \cite{IQN}, the incremental quasi-Newton method (IQN) is proposed to minimize the objective function written in  a sum of large amounts of strongly convex functions. It alleviates the high computational cost at each iteration. The main ingredients are as follows. In lieu of random selection of an individual function, incremental methods choose this individual function in a cyclic routine. Thus, it leads to efficient implementation of both the BFGS and iterate updates. The aggregated gradients of all functions are successful in reducing the noise of gradient approximation. Moreover, it satisfies the  Dennis-Mor\'{e} condition. This indicates that IQN method yields local superlinear convergence rate.

Furthermore, the quality of the curvature estimate may be difficult to control in stochastic regime. To alleviate it,  \cite{RHBYRD} has investigated an efficient subsampled Hessian-vector product  to estimate the curvature information based on the limited memory BFGS (LBFGS). This method is applied in strongly convex optimization and can avoid doubly evaluating gradients. In  \cite{RHBYRD2}, the subsampled Hessian matrix scheme is adopted in matrix-vector product form, and the conjugate gradient method is further applied to obtain the search direction. Moreover, the subsampled Hessian matrix is also used as the initial Hessian approximation matrix in the LBFGS method. This is because the traditional choice contains little curvature information about the problem. In  \cite{SBBFGS}, the subsampled Hessian matrix has been adopted to formulate the stochastic block BFGS scheme. The main ingredient is left-multiplying the inverse equation by a randomly generated matrix with few columns. Hence, the computational cost is substantially reduced. In \cite{NIPS2013_4937}, the stochastic variance reduced gradient (SVRG) strategy  has been employed to reduce the variance of the stochastic gradient.  

It should be noted that the above discussed second-order methods have been proposed for solving convex problems. They cannot be directly applied to nonconvex problems. Moreover, tackling non-convexity and ill-conditioning are two major challenges in stochastic nonconvex optimization problems. To this end, damped BFGS  \cite{SQN} and regularized BFGS \cite{MokhtariAryan2014RRSB} have been proposed to deal with the non-convexity and ill-conditioning of the stochastic optimization problem, respectively. In stochastic BFGS methods, the Hessian approximation matrices are ensured to be positive definite in strongly convex optimization problems  \cite{nocedal2006numerical}. However, it is not the case for nonconvex objective functions. In \cite{SQN}, a stochastic damped BFGS based on  \cite{Powell1978} is proposed to address this issue. However, the BFGS update may still be ill-conditioned if there are insufficient samples. Moreover, the convergence may be significantly affected if the BFGS matrix is close to singularity or even singular. In  \cite{MokhtariAryan2014RRSB}, a regularized stochastic BFGS (RES) method is proposed to improve the numerical condition mentioned above. However, if the problem is nonconvex, the BFGS update may become non-positive definite and hence a descent step may not be guaranteed. Moreover, directly combining the damped scheme  \cite{MokhtariAryan2014RRSB} and this regularized formulation may still not be  able to guarantee positive definiteness of the BFGS update and a descent step. 
To this end, we propose in this paper a novel stochastic quasi-Newton method, called Sd-REG-LBFGS method, to address the above problems. Our main contributions are as follows:

\begin{itemize}
	\item New damped BFGS scheme: We propose a new stochastic regularized damped BFGS method containing a novel damped parameter and a new gradient difference scheme. The proposed scheme guarantees positive definiteness of the BFGS update and improves the numerical condition of the optimization problem. 
	\item Choice of Regularization Parameters: The choice of the regularization parameters for the new regularized gradient difference and damped parameter schemes is crucial to ensure positive definiteness of the BFGS update. We proved that if the chosen regularization parameters satisfy a certain condition (Lemma 1) we have derived, then positive definiteness is guaranteed for the proposed approach.
	\item Convergence Analysis: The convergence property of the proposed method is thoroughly analyzed. In particular, we show that the norm of the updated Hessian approximation matrix is uniformly bounded (see Lemmas 2 and 3), which is a necessary condition for convergence. Furthermore, we showed that with a specified step size, the iteration number $N$ required to reach a norm of gradient of  $\frac{1}{N}\sum_{k=0}^{N-1}\mathbb{E}(\norm{\nabla f(x_k)}^2)<\epsilon$ is  at most $O(\epsilon^{-\frac{1}{1-\upsilon}})$, for $0.5<\upsilon<1$. All the above convergence results are independent of the convexity assumption. Thus, our proposed method can be applied to nonconvex problems.
	
\end{itemize}

For numerical study, the proposed approach is evaluated using a logistic regression,  a Bayesian logistic regression and a nonconvex relaxed soft margin support vector machine (SVM)\footnote{  Due to page limitation, the simulation results for the nonconvex relaxed soft margin SVM is omitted here and interested readers are referred to Section V of the supplementary material.}. Experimental results using a synthetic dataset and several real datasets  \cite{scene,Elec,wifi1,wifi2,IONO,BNA}  show that the proposed regularized damped stochastic BFGS method performs better than the conventional damped stochastic BFGS and other algorithms in terms of classification accuracy (ACC) and norm of gradient (NOG), which suggest it converges closer to the stationary point. Moreover, the sensitivity of the proposed algorithm on various algorithmic parameters and the complexity of the proposed algorithm are also studied. Due to page limitation, it is omitted here and interested readers are referred to Sections III and IV of the supplementary material for details.

The rest of the paper is organized as follows: Section II reviews the general formulation of the SQN framework. In Section III, we provide the detail derivation of our proposed algorithm, including the uniform bound on the norm of LBFGS matrix and the convergence results. In Sections IV and V, the effectiveness of the proposed Sd-REG-LBFGS algorithm is demonstrated through solving several machine learning problems, and the numerical experiments are conducted to evaluate the performance of the proposed algorithm with a comparison with conventional algorithms. The conclusion is provided in Section VI.

\textit{Mathematical Notation}: we use $\norm{a}$ to denote the Euclidean norm of vector $a$ and $\norm{A}$ to denote the matrix norm of a matrix $A$. The trace operator of $A$ is written as $\text{Tr}(A)$ and the determinant as $\text{det}A$. The operator $\mathbb{E}_{\Xi}(\cdot)$ stands for the expectation taken with respect to random variable $\Xi$. $A\succeq B$ indicates the matrix $A-B$ is positive semidefinite. The identity matrix with appropriate dimension is signified as $I$. 

\section{Problem Formulation}

Consider the following general optimization problem in expectation form:
\begin{equation}\label{optprob1}
\underset{x\in{\mathbb{R}^n}}{\text{min}}\quad f(x):=\mathbb{E}[F(x,\Xi)],
\end{equation} 
where $\Xi\in{\mathbb{R}^d}$ denotes a random variable, and $F: \mathbb{R}^n\times\mathbb{R}^d\rightarrow{\mathbb{R}}$ is possibly a nonconvex random function. In many applications, the expectation in (\ref{optprob1}) is intractable, or the value and gradient of $f$ are not easily obtained. For example, in machine learning problems, the random variables may contain the input features $Y$ and the class labels $Z$, i.e.  $\Xi=(Y,Z)$, which may follow some unknown distribution $P$, in which inferences are to be made. The training set is assumed to be a collection of independent and identically distributed (i.i.d.) samples $\xi_i=(y_i,z_i)$ with $i=1,\dots,N$, distributed according to $P$  via certain observations. The expectation of $F(x,\Xi)$ in (\ref{optprob1}) can be approximated by the following empirical average  $\bar{f}(x)=1/N\sum_{i=1}^{N} F (x,\xi_i)$, where $F (x,\xi_i)$ is the empirical loss function corresponding to the same $i$-th sample $\xi_i$. For a  large-scale problem where  $N$ is large, this exact empirical gradient may require expensive evaluation of  $F (x,\xi_i)$  for all the samples. In general, stochastic optimization can also be applied to problems where one might be able to access values of the objective function and its gradient from some physical sensor devices in physical simulations. The measured results may be noisy and depend on the unknown $\xi_n$ every time we attempt to measure $F(x,\xi)$ or its gradient.

In this paper, we mainly focus in machine learning problems mentioned above. Moreover, the stochastic gradient, denoted as $g(x,\Xi)$ is an unbiased estimator of $\nabla{f(x)}$, i.e., $\mathbb{E}_{\Xi}[g(x,\Xi)]=\nabla{f(x)}$, where the expectation is taken with respect to $\Xi$.   We assume that we can access the gradient via explicit evaluation from the training data (or some physical sensor devices in physical simulations for general stochastic optimization).  In addition, we assume that  $f$ is continuously differentiable and  the gradient of $f$ is Lipschitz continuous:
\begin{equation}\label{gradientLipschitz}
\norm{\nabla{f(x)}-\nabla{f(y)}} \leq{L_f\norm{x-y}},
\end{equation}
with Lipschitz constant $L_f>0$. 

In classical deterministic quasi-Newton methods, at iteration $k$, the update of current iterate is given by: 
\begin{equation}\label{QN}
x_{k+1}=x_k-\eta_kB^{-1}_k\nabla{f(x_k)},
\end{equation}
where $B_k$ is an approximation to the Hessians of the objective
function $\nabla^2f(x_k)$, since evaluating $\nabla^2f(x_k)$ is computationally intensive. Various Hessian approximation methods have been proposed which include, e.g., Broyden, Fletcher, Goldfarb, and Shanno (BFGS);  Davidon, Fletcher, and Powell (DFP) and symmetric rank-1 (SR1) updates. In this paper, we mainly focus on the following BFGS update  as it is one of the most popular quasi-Newton algorithms:
\begin{equation}\label{BFGS}
B_{k+1}=B_{k}+\frac{y_{k}y_{k}^T}{s^T_{k}y_{k}}-\frac{B_{k}s_{k}s^T_{k}B_{k}}{s^T_{k}B_{k}s_{k}},
\end{equation}
where the correction pairs  are  $s_{k}=x_{k+1}-x_{k}$ and $y_{k}=\nabla{f(x_{k+1})}-\nabla{f(x_{k})}$ respectively.  It can be shown that (\ref{BFGS}) satisfies the secant equation, i.e., $B_ks_{k+1}=y_{k}$. To show that the resultant matrix is positive definite, one can rewrite  (\ref{BFGS}) by letting  $s=s_{k}, y=y_{k}, B=B_{k}, B_{k+1}=B^+$ for notational convenience, which yields:
\begin{equation}\label{66}
B^+=\frac{yy^T}{s^Ty}+B^{\frac{1}{2}}\left(I-\frac{B^{\frac{1}{2}}ss^TB^{\frac{1}{2}}}{s^TB^{\frac{1}{2}}B^{\frac{1}{2}}s}\right)B^{\frac{1}{2}}.
\end{equation}   
Moreover, it can be shown by induction that with the condition $s^T_{k}y_{k}>0$, and an initial positive definite Hessian approximation $B_0\succ{0}$, $B_k$ is updated recursively and remains positive definite in subsequent iterations. In fact, the condition $s^T_{k}y_{k}>0$ to preserve the positive definiteness of the Hessian approximation update via (\ref{BFGS}) is always satisfied for strongly or strictly convex objective functions. This is due to the monotonic gradient mapping property \cite{CVX}. To be specific, if the objective function $f$ is strongly or strictly convex, for any $x,y\in{\mathbb{R}^n}$, $(\nabla{f(x)}-\nabla{f(y)})^T(x-y)>0$. Hence, by letting $x=x_{k+1}$ and $y=x_{k}$, we can see that the condition $s^T_{k}y_{k}>0$ is satisfied.

To migrate the classical quasi-Newton method to the stochastic regime, the main ingredient is to adopt the stochastic approximation for the exact gradient, which forms the general framework of the SQN method. More precisely, at iteration  $k$, we subsample a mini-batch $m_k$ of data so as to compute  the stochastic  gradient evaluated at the current solution  $x_k$, which we shall refer to as $\nabla{F}(x_k,\xi_{k,i})$  with $i=1,\dots,m_k$. The SA based on this mini-batch estimate can be obtained by the following ensemble average of $\nabla{F}(x_k,\xi_{k,i})$ with $i=1,\dots,m_k$: $\bar{g}(x,\xi_k)=\frac{1}{m_k}\sum_{i=1}^{m_k}\nabla{F}(x_k,\xi_{k,i})$. By combining (\ref{QN}) and (\ref{BFGS}), one gets the desired SQN iterate as follows: 
\begin{equation}
x_{k+1}=x_k-\eta_k{B^{-1}_k}\bar{g}(x_k,\xi_k),
\end{equation}  
where the following stochastic gradient difference is employed in BFGS update (\ref{BFGS}):
\begin{equation}\label{SGDif}
y_{k}=\frac{1}{m_{k}}\sum_{i=1}^{m_{k}}\nabla{F}(x_{k+1},\xi_{k,i})-\nabla{F}(x_{k},\xi_{k,i}).
\end{equation}

\begin{remark}
	It should be noted from the first term in (\ref{SGDif}) that the  gradient of $F$ at $x_{k+1}$ is generated using the same  subsampling process conducted at current iteration. This implies that at each iteration, the stochastic gradient is evaluated twice. There are two advantages: i). For strongly convex function $F(\cdot)$, using (\ref{SGDif}) guarantees the condition $s^T_ky_k>0$. Moreover, we suggest to adopt the first-order Taylor approximation to reduce the computational complexity, i.e., $y_k\approx1/m_k\sum_{i=1}^{m_k}\nabla^2F(x_k,\xi_{k,i})s_k$, where  $\nabla^2F(x_k,\xi_{k,i})s_k$ is a product of the matrix and the vector, which can be  obtained with low complexity  \cite{RHBYRD}; ii). It ensures that the BFGS Hessian approximations are uniformly bounded below and above.
\end{remark}

\section{The Proposed Algorithm}
\subsection{The Proposed Damped SQN Method}
In nonconvex optimization problems, the positivity condition  $s^T_ky_k>0$ of the correction pairs may not be maintained. This may lead to non-positive  definite BFGS matrix. To remedy this problem,  \cite{Powell1978} has proposed a damped QN method to preserve the positive definiteness of  BFGS matrix in nonconvex optimization. Here, we shall extend it to stochastic regime. Specifically, $y_k$ is modified to $\bar{y}_{k}:=\theta_{k}y_{k}+(1-\theta_{k})B_{k}s_{k}$ (thus $y_k$   in (\ref{SGDif}) will be modified), where $\theta_{k}$ is the damped parameter satisfying:
\begin{equation}\label{damp1}
\theta_{k}=\left\{
\begin{aligned}
&\frac{0.8s^T_{k}B_{k}s_{k}}{s^T_{k}B_{k}s_{k}-s^T_{k}y_{k}}, && \text{if } s^T_{k}y_{k}\leq 0.2s^T_{k}B_{k}s_{k},\\
&\quad\quad\quad 1, && \text{otherwise}.
\end{aligned}
\right.
\end{equation}

It can be easily verified that $B_k\succ{0}$ and $0<\theta_{k}\leq{1}$ with an initial positive definite Hessian approximation $B_0\succ 0$. Note that when $\theta_{k}=1$, which is often the case in practice, the BFGS matrix update reduces to the classical formula in (\ref{BFGS}). For other values of $\theta_k$, such modification prevents the determinant of $B_{k+1}$ from being less than 0.2 of the determinant of $B_k$ \cite{Powell1978}. In addition, since:
\begin{equation}\label{damp}
s^T_k\bar{y}_{k}=\left\{
\begin{aligned}
&0.2s^T_kB_ks_k, && \text{if } s^T_{k}y_{k}\leq 0.2s^T_{k}B_{k}s_{k},\\
&s^T_ky_k, && \text{otherwise},
\end{aligned}
\right.
\end{equation}     
it implies that if $B_k\succ{0}$, then $s^T_k\bar{y}_{k}\geq{0.2s^T_kB_ks_k}>0$, and the damped quasi-Newton method ensures the  positive definiteness of the BFGS update $B_{k+1}$. 

For nonconvex optimization problems, even the stochastic damped BFGS method guarantees all the subsequent $B_{k+1}$ obtained via  (\ref{BFGS}) be positive definite, it is possible for the smallest eigenvalue of $B_{k+1}$ to be arbitrarily close to zero, and hence, the Hessian approximation matrix $B_k$ will be nearly singular \cite{MokhtariAryan2014RRSB}. To remedy the problem, we shall propose a generalized RES scheme for nonconvex optimization using novel damped QN method. We shall first introduce briefly the regularized stochastic quasi-Newton method (RES) for strongly convex optimization problems in \cite{MokhtariAryan2014RRSB}. Then, the proposed generalized RES scheme will be described. 

Recall $B_{k+1}$ in (\ref{BFGS}) is obtained by solving the following semidefinite programming problems:
\begin{equation}\label{BFGSopt}
\begin{aligned}
& \underset{Z}{\text{min}}
& & \text{Tr}[B^{-1}_kZ]-\text{logdet}[B^{-1}_kZ]-n \\
& \: \text{s.t.}
& & Zs_k=y_k, Z\succeq 0,
\end{aligned}
\end{equation}
where the optimal solution to (\ref{BFGSopt}) is $Z^*=B_{k+1}$, obtained by nulling the gradient of the Lagrangian duality function $\varphi(Z(\nu),\nu)=\text{inf}_{Z\succeq{0}}\;\mathcal{L}(Z,\nu)$ with respect to $\nu$, in which $\mathcal{L}(Z,\nu)=\text{Tr}[B^{-1}_kZ]-\text{logdet}[B^{-1}_kZ]-n+\nu^T(Zs_k-y_k)$. A simple interpretation to (\ref{BFGSopt}) is to minimize the Gaussian differential entropy between the Gaussian distributions $\mathcal{N}(0,B_k)$ and $\mathcal{N}(0,Z)$ with the constraint of the secant equation and positive semidefinite solution. For the RES strategy, the following modification of the optimization problem (\ref{BFGSopt}) is solved:
\begin{equation}\label{BFGSopt2}
\begin{aligned}
& \underset{Z}{\text{min}}
& & \text{Tr}[B^{-1}_k(Z-\gamma I)]-\text{logdet}[B^{-1}_k(Z-\gamma I)]-n \\
&\: \text{s.t.}
& & Zs_k=y_k, Z\succeq 0.
\end{aligned}
\end{equation}

By setting $\tilde{Z}=Z-\gamma I$ and $\tilde{y}_k=y_k-\gamma s_k$,  
the following regularized BFGS update is obtained by using the related Lagrangian duality function:
\begin{equation}\label{res}
B_{k+1}=B_{k}+\frac{\tilde{y}_{k}\tilde{y}^T_{k}}{s^T_{k}\tilde{y}_{k}}-\frac{B_{k}s_{k}s^T_{k}B_{k}}{s^T_{k}B_{k}s_{k}}+\gamma I.
\end{equation}

Under the condition $s^T_{k}\tilde{y}_{k}>0$ with an initial positive semidefinite $B_0\succeq 0$, the subsequent Hessian approximations will have the smallest eigenvalue exceeding a given desired level  $\gamma$. Comparing (\ref{BFGS}) and (\ref{res}), one can see that not only is $y_k$ being modified to $\tilde{y}_k$, an additional regularization term $\gamma I$ is also introduced to avoid possible ill-conditioning.

However, it can be verified that RES cannot be adopted to the damped QN mehtod for nonconvex optimization problems by simply applying (\ref{damp1}) to modify $y_k$ in (\ref{res}). We briefly illustrate this below. Consider $\bar{y}_k$, which is the modified version of $y_k$ by employing (\ref{damp1}). It follows that $s_k^T\tilde{y}_k$ can be calculated as follows:
\begin{equation}\label{illus}
s^T_k\tilde{y}_{k}=\left\{
\begin{aligned}
&0.2s^T_kB_ks_k-\gamma s^T_ks_k, && \text{if } s^T_{k}y_{k}\leq 0.2s^T_{k}B_{k}s_{k},\\
&s^T_ky_k-\gamma s^T_ks_k, && \text{otherwise},
\end{aligned}
\right.
\end{equation}     
Hence, the positivity of $s_k^T\tilde{y}_k$ cannot be guaranteed. Moreover, 
even in strongly convex functions $F(\cdot)$ with convexity parameter $\underline{m}$ (i.e., $\nabla^2F\succeq\underline{m}I$), if the given level $\gamma$ is chosen to be greater than $\underline{m}$, which results in $s^T_k\tilde{y}_k<0$, $B_{k+1}$ can still be near singular or negative positive. 

To remedy the problem, we now propose a novel damped SQN method. To start with, the following stochastic gradient difference $\hat{y}_k$ is proposed to modify $y_k$:
\begin{equation}\label{sgd2}
\hat{y}_k=\bar{\theta}_{k}y_k+(1-\bar{\theta}_{k})(B_k+\delta I)s_k,
\end{equation}
where $\delta$ is a given positive constant that satisfies specific condition (see Lemma 1). Furthermore, we propose to update the damped parameter as follows: 
\begin{equation}\label{damp2}
\bar{\theta}_{k}=\left\{
\begin{aligned}
&\frac{0.8s^T_{k}(B_{k}+\delta I)s_{k}-\gamma s^T_ks_k}{s^T_{k}(B_{k}+\delta I)s_{k}-s^T_{k}y_{k}}, &&\begin{aligned}\text{if }
 &s^T_{k}y_{k}\leq 0.2s^T_{k}(B_{k}\\
 &+\delta I)s_{k}+\gamma s^T_ks_k,
\end{aligned}\\
&\quad\quad\quad\quad\quad 1, && \text{otherwise}.
\end{aligned}
\right.
\end{equation}
Substituting $\tilde{\hat{y}}_k:=\hat{y}-\gamma s_k$  into (\ref{res}) with the parameter $\bar{\theta}_{k}$ defined in (\ref{damp2})  yields our proposed Hessian approximation updating scheme:
\begin{equation}\label{res2}
B_{k+1}=B_{k}+\frac{\tilde{\hat{y}}_{k}\tilde{\hat{y}}^T_{k}}{s^T_{k}\tilde{\hat{y}}_{k}}-\frac{B_{k}s_{k}s^T_{k}B_{k}}{s^T_{k}B_{k}s_{k}}+\gamma I.
\end{equation}
The following lemma shows that by recursively updating $B_k$ via (\ref{res2}), our proposed method  maintains the positive definiteness of the Hessian approximation matrix at each iteration. 

\textit{Lemma 1}. For $\hat{y}_k$ defined in (\ref{sgd2}) and $\delta$ is chosen to satisfy $0.8\delta\geq{\gamma}$, then $0<\bar{\theta}_{k}\leq{1}$ and $s^T_k\tilde{\hat{y}}_k\geq{0.2s^T_k(B_k+\delta I)s_k}$. Moreover, if $B_{k}\succ{0}$, then $B_{k+1}$ generated by the proposed damped BFGS update (\ref{res2}) are positive definite with the smallest eigenvalue exceeding the given desired level $\gamma$.

\textit{Proof.} Note from (\ref{damp2}) that, if $s^T_{k}y_{k}\leq 0.2s^T_{k}(B_{k}+\delta I)s_{k}+\gamma s^T_ks_k$, then $\bar{\theta}_k=1$; for $s^T_{k}y_{k}> 0.2s^T_{k}(B_{k}+\delta I)s_{k}+\gamma s^T_ks_k$, by substituting the inequality into $\bar{\theta}_k$, we get the following inequality: 
\begin{equation}
\begin{aligned}
\bar{\theta}_k&=\frac{0.8s^T_{k}(B_{k}+\delta I)s_{k}-\gamma s^T_ks_k}{s^T_{k}(B_{k}+\delta I)s_{k}-s^T_{k}y_{k}}\\
&\leq \frac{0.8s^T_{k}(B_{k}+\delta I)s_{k}-\gamma s^T_ks_k}{s^T_{k}(B_{k}+\delta I)s_{k}-[0.2s^T_{k}(B_{k}+\delta I)s_{k}+\gamma s^T_ks_k]}=1.
\end{aligned}
\end{equation}
Moreover, the numerator of (\ref{damp2}) satisfies  $0.8s^T_{k}B_{k}s_{k}+(0.8\delta-\gamma) s^T_{k}s_{k}\geq 0.8s^T_{k}B_{k}s_{k} >0$ with the conditions $0.8\delta\geq\gamma$ and $B_k\succ 0$. Similarly from the denominator in (\ref{damp2}), we have:
\begin{equation}
s^T_{k}(B_{k}+\delta I)s_{k}-s^T_{k}y_{k}\geq 0.8s^T_{k}(B_{k}+\delta I)s_{k}-\gamma s^T_ks_k>0.
\end{equation}    
Subsequently, both the numerator and denominator of (\ref{damp2}) are positive and its maximum value is one, i.e.,  $0<\hat{\theta}_k\leq 1$. Moreover, from (\ref{sgd2}) and (\ref{damp2}), $s^T_k\tilde{\hat{y}}_k$ can be calculated as follows:
\begin{equation}\label{eq22}
\begin{aligned}
s^T_k\tilde{\hat{y}}_{k}&=s^T_{k}(B_{k}+\delta I)s_{k}-\gamma s^T_ks_k-\bar{\theta}_k[s^T_{k}(B_{k}+\delta I)s_{k}-s^T_{k}y_{k}]\\
&=\left\{
\begin{aligned}
&0.2s^T_k(B_k+\delta I)s_k, &&\begin{aligned} \text{if }
s^T_{k}y_{k}&\leq 0.2s^T_{k}(B_{k}+\delta I)s_{k}\\
&\quad\quad+\gamma s^T_ks_k,
\end{aligned} \\
&s^T_ky_k-\gamma s^T_ks_k, && \text{otherwise}.
\end{aligned}
\right.
\end{aligned}
\end{equation}     
From (\ref{eq22}), we can see that $s^T_k\tilde{\hat{y}}_{k}\geq 0.2s^T_k(B_k+\delta I)s_k$. Therefore, if $B_k$ is positive definite, it follows that $s^T_k\tilde{\hat{y}}_{k}>0$. Consequently, as in (\ref{66}), the first three terms in the right hand side of the proposed BFGS update scheme (\ref{res2}) is a positive definite matrix.          
\begin{remark}
	From the inequality $s^T_k\tilde{\hat{y}}_{k}\geq 0.2s^T_k(B_k+\delta I)s_k$, we further have $s^T_k\tilde{\hat{y}}_{k}\geq 0.2[\lambda(B_k)_{min}+\delta]s^T_ks_k$, where $\lambda(B_k)_{min}$ is the smallest eigenvalue of $B_k$. Next, we shall extend the proposed BFGS update to a limited memory version.
\end{remark} 
\subsection{The Proposed Algorithms for Limited Memory }

The limited-memory quasi-Newton method \cite{Liu1989}, which approximates the Hessian approximation from a limited number of vectors attained from recent iterations, is useful in large scale applications to reduce the large memory storage of the Hessian approximation matrices. As this method requires  modest storage and possesses good convergence speed, it is generally considered to be superior to the steepest descent method for deterministic optimization \cite{RHBYRD}. Interested readers are referred to \cite{nocedal2006numerical} for more information.  In recent years, stochastic limited-memory BFGS (L-BFGS) methods  have been studied for strongly convex optimization problems \cite{pmlr-v51-moritz16}\cite{Liu1989}\cite{RHBYRD}. In this subsection, we propose a stochastic  damped and regularized L-BFGS (Sd-REG-LBFGS) method for nonconvex optimization problems. 

For robustness in implementation and to amortize the cost, one of the strategies is to update the BFGS Hessian approximation at spaced intervals using the average of the iterate points instead of at each iteration \cite{RHBYRD}. Motivated by this strategy, we compute the correction pairs $\{s_j,y_j\}$ based on the average of the iterates in the specified interval. The BFGS Hessian approximations are subsequently calculated. In particular, all the modifications are based on our proposed damped BFGS method in (\ref{sgd2})-(\ref{res2}).  Specifically, we assume that the length of the aforementioned interval of iterations is $L$. Suppose we have a memory with size $M$. It stores the sequence of correction pairs $\{s_j,y_j\}$ for $j=t-(M-1)-1,\dots,t-1$, where $t:=\frac{k+1}{L}$ and the iteration $k$ satisfies $(k+1)\;\text{mod}\;L=0$ and $k\geq M(L-1)-1$. We further define $s_j$ as the difference of two average iterates with respect to the two most recent disjoint intervals, i.e.,:
\begin{equation}\label{displacement}
s_{j}=\bar{x}_{j+1}-\bar{x}_j, \; \text{where}\; \bar{x}_{j}=\left\{
\begin{aligned}
&\frac{1}{L}\sum_{k=(j-1)L}^{jL-1}x_k, && \text{if \;}j\geq 1,\\
&\quad x_0, && \text{if\;\;} j=0.
\end{aligned}
\right. 
\end{equation}
Subsequently, the gradient difference is evaluated at $\bar{x}_{j+1}$ and $\bar{x}_j$ as follows:
\begin{equation}\label{sgd3}
y_j=\frac{1}{m_j}\sum_{l=1}^{m_j}\nabla{F(\bar{x}_{j+1},\xi_{j,l})}-\nabla{F(\bar{x}_{j},\xi_{j,l})}.
\end{equation} 

Recall that we only update BFGS matrix at the end of each interval, to reduce the memory of storing $B_t$, we can further approximate it using the L-BFGS method, where a  sequence of correction pairs in (\ref{displacement}) and (\ref{sgd3}) are stored. Based on the stochastic damped and regularized BFGS method proposed in (\ref{res2}), we define a new vector $\tilde{y}_j:=\hat{\theta}_jy_j+(1-\hat{\theta}_j)(\hat{B}^{(0)}_{j+1}+\delta I)s_j-\gamma s_j$, with $\hat{\theta}_j$ given by:
\begin{equation}\label{damp3}
\hat{\theta}_{j}=\left\{
\begin{aligned}
&\frac{0.8s^T_{j}(\hat{B}^{(0)}_{j+1}+\delta I)s_{j}-\gamma s^T_js_j}{s^T_{j}(\hat{B}^{(0)}_{j+1}+\delta I)s_{j}-s^T_{j}y_{j}}, && \begin{aligned}
&\text{if } s^T_{j}y_{j}\leq \gamma s^T_js_j+\\
&0.2s^T_{j}(\hat{B}^{(0)}_{j+1}+\delta I)s_{j},
\end{aligned}\\
&\quad\quad\quad\quad\quad 1, && \text{otherwise},
\end{aligned}
\right.
\end{equation}
where $\hat{B}^{(0)}_{j+1}$ is an initial estimate of the Hessian matrix and a typical value of $\hat{B}^{(0)}_{j+1}$ in standard L-BFGS is $\frac{y^T_{j}y_{j}}{s^T_{j}y_{j}}I$. As the denominator $s^T_{j}y_{j}$ may not be positive for nonconvex problems, we propose the following initial value of $\hat{B}^{(0)}_{j+1}$:
\begin{equation}\label{B}
\hat{B}^{(0)}_{j+1}=\tau_{j+1}I,\;\text{where}\;\tau_{j+1}=\text{max}\left\{\frac{y^T_{j}y_{j}}{s^T_{j}y_{j}}+\gamma,\beta\right\} ,
\end{equation}
where $\beta$ is a given positive constant and is also the lower bound on $\tau_{j}$, i.e., $\tau_{j}>\beta$. Therefore, at the end of the $t$-th interval, we define the Sd-REG-LBFGS formula from the past correction pairs $(s_j,\tilde{y}_j)$ via the following inner iterations:
\begin{equation}\label{res3}
\hat{B}^{(i+1)}_{t}=\hat{B}^{(i)}_{t}+\frac{\tilde{y}_{j}\tilde{y}^T_{j}}{s^T_{j}\tilde{y}_{j}}-\frac{\hat{B}^{(i)}_{t}s_{j}s^T_{j}\hat{B}^{(i)}_{t}}{s^T_{j}\hat{B}^{(i)}_{t}s_{j}}+\gamma I
\end{equation}   
for $i=0,\dots,M-1$ and $j=t-(M-1)+i-1$. It follows from \textit{Lemma 1}  that $s^T_j\tilde{y}_j\geq 0.2s^T_j(B^{(0)}_{j+1}+\delta I)s_j$. Therefore, starting with the positive definite matrix $\hat{B}^{(0)}_{t}$ given in (\ref{B})  and a constant $\delta$ satisfying $0.8\delta>\gamma$, the positive definite matrix $\hat{B}_{t}=\hat{B}^{(M)}_{t}\succ \gamma I$ can be updated by the inner iteration of the  proposed Sd-REG-LBFGS formula in (\ref{res3}). Furthermore, as the gradient is stochastic and the exact evaluation of the objective function is expensive at each iteration, the Wolfe condition based on the incomplete stochastic gradient may lead to premature condition for convergence or oscillation and prevent the algorithm further progressing. Therefore, we choose the step size to satisfy the well-known condition \cite{RobbinsHerbert1951ASAM} for the step size choice in stochastic optimization, namely: 
\begin{equation}
	\sum_{k=1}^{\infty}\eta_k=\infty,\quad\sum_{k=1}^{\infty}\eta^2_k<\infty.
\end{equation}
A popular choice is $\eta_k=\frac{r}{k}$, for $r>0$ \cite{MokhtariAryan2014RRSB,SQN,RHBYRD}.
The proposed Sd-REG-LBFGS algorithm is summarized in Algorithm 1.

\begin{algorithm}[t]
	\renewcommand{\algorithmicrequire}{\textbf{Input:}}
	\renewcommand{\algorithmicensure}{\textbf{Output:}}
	\caption{Sd-REG-LBFGS} 
	\label{alg1}
	\begin{algorithmic}[1]
		\REQUIRE initial optimization variable $x_0$, memory size $M$, interval length $L$, step length  \hspace*{0.02in} $\eta_k$ and gradient sample batch size $m_k$, choose the constant $\delta$ and $\gamma$ satisfying $0.8\delta>\gamma$ 
		\STATE Set $t=0$ and generate $m_{0}$ samples $\{\xi_{0,l}\}^{m_{0}}_{l=1}$ 
		\FOR{$k=0,1,\dots$}
		\STATE Randomly choose $m_k$ samples $\xi_k=\{\xi_{k,1},\cdots,\xi_{k,m_k}\}$
		\STATE Calculate stochastic gradient $\bar{g}(x_k,\xi_k)=\frac{1}{m_k}\sum_{l=1}^{m_k}\nabla F(x_k,\xi_{k,l})$,
		
		\IF{$t<2$} 
		\STATE$x_{k+1}=x_k-\eta_k\bar{g}(x_k,\xi_k)$
		\ELSE
		\STATE$x_{k+1}=x_k-\eta_k\hat{B}^{-1}_{t}\cdot\bar{g}(x_k,\xi_k)$ 
		\ENDIF
		\IF{$(k+1)\;\text{mod}\;L=0$}

		\STATE  Calculate and store the correction pairs:   $s_{t}$ and $y_{t}$ according to (\ref{displacement}) and (\ref{sgd3}) respectively
		\STATE Set $t=t+1$
		\STATE  Generate $m_{t}$ samples $\{\xi_{t,l}\}^{m_{t}}_{l=1}$
		
		\IF{$t>1$}
		\STATE Set $\tilde{M}=\text{min}\{t,M\}$, draw the sequence of correction pairs $\{ s_j,y_j \}^{t-1}_{j=t-\tilde{M}}$ from the memory.
		\STATE Set the initial matrix $\hat{B}^{(0)}_{t}=\tau_{t}I,\;\text{where}\;\tau_{t}=\text{max}\left\{\frac{y^T_{t-1}y_{t-1}}{s^T_{t-1}y_{t-1}}+\gamma,\beta\right\} $
		\FOR{$i=0,\dots,\tilde{M}-1$}
		\STATE Set $j=t-\tilde{M}+i$ and apply Sd-REG-LBFGS formula according to (\ref{res3})
		\ENDFOR
		\STATE Set $\hat{B}_t=\hat{B}^{(\tilde{M})}_t$.
		\ENDIF

		\ENDIF
		
		\ENDFOR
		
	\end{algorithmic}
\end{algorithm}

\subsection{Convergence Result}
For the convergence result of our proposed algorithm, one significant condition is that the norm of the resulting $\hat{B}^{(i+1)}_{t}$ from (\ref{res3}) is uniformly bounded above, and uniformly bounded below from zero. Moreover, the following assumption is useful for the derivation of the upper and lower bound:

\textit{Assumption 1}\cite{SQN}. The random function $F(x,\Xi)$ is twice continuously differentiable, where the second-order derivative with respect to $x$ is denoted as $\nabla^2F(x,\Xi)$. Moreover, there exists a positive constant $\rho$ such that $\norm{\nabla^2F(x,\Xi)}\leq \rho$.

Note that the above assumption implies that $-\rho I\prec\nabla^2F(x,\Xi)\prec \rho I$, rather than the strong convexity assumption $0\prec\underline{\rho} I\prec\nabla^2F(x,\Xi)\prec \bar{\rho} I$ in \cite{RHBYRD}\cite{MokhtariAryan2014RRSB}. The following lemma shows that the norm of the matrix $\hat{B}^{\tilde{M}}_t$ generated by the Sd-REG-LBFGS formula (\ref{res3}) is uniformly bounded above.

\textit{Lemma 2}.  Given the positive definite matrix $\hat{B}^{(0)}_t$ defined by (\ref{B}), suppose $\hat{B}^{(i+1)}_t$ is updated through L-BFGS computation step in the $t$-th interval of Algorithm 1, then with Assumption 1, the norm of $\hat{B}^{(\tilde{M})}_t$ is bounded above, i.e.,
\begin{equation}\label{B_aboveBound}
\norm{\hat{B}^{(\tilde{M})}_t}\leq Q_U,
\end{equation} 
where $Q_U=\beta+\rho+\gamma+\tilde{M}(Q+5\rho+\gamma) $, $\tilde{M}=\text{min}\{t,M\}$ and $Q$ is defined as follows:
\begin{equation}\label{Q}
Q=\text{max}\left\{\begin{aligned}\frac{5(\rho+\gamma)^2}{\beta+\delta}+5(\beta+\delta), &\frac{5(\rho+\gamma)^2}{\beta+\rho+\gamma+\delta}+\\&\quad 5(\beta+\rho+\gamma+\delta)\end{aligned} \right\}.
\end{equation}

\textit{Proof}. Recall from the Sd-REG-LBFGS formula that according to \textit{Lemma 1}, each generated matrix satisfies $\hat{B}^{(i+1)}_t\succ \gamma I$. Note from the third term on the right hand side in (\ref{res3}) that  the matrix term $\frac{\hat{B}^{(i)}_{t}s_{j}s^T_{j}\hat{B}^{(i)}_{t}}{s^T_{j}\hat{B}^{(i)}_{t}s_{j}}$ is positive definite. Therefore, we have:
\begin{equation}\label{proofeq1}
\hat{B}^{(i+1)}_{t}\preceq\hat{B}^{(i)}_{t}+\frac{\tilde{y}_{j}\tilde{y}^T_{j}}{s^T_{j}\tilde{y}_{j}}+\gamma I.
\end{equation}  
Taking matrix norm on both sides and using triangle inequality of norm leads to:
\begin{equation}\label{proofeq2}
\begin{aligned}
\norm{\hat{B}^{(i+1)}_{t}}&\leq\norm{\hat{B}^{(i)}_{t}+\frac{\tilde{y}_{j}\tilde{y}^T_{j}}{s^T_{j}\tilde{y}_{j}}+\gamma I}\leq\norm{\hat{B}^{(i)}_{t}}+\norm{\frac{\tilde{y}_{j}\tilde{y}^T_{j}}{s^T_{j}\tilde{y}_{j}}}+\gamma\\
&=\norm{\hat{B}^{(i)}_{t}}+\frac{\tilde{y}^T_{j}\tilde{y}_{j}}{s^T_{j}\tilde{y}_{j}}+\gamma,
\end{aligned}
\end{equation}
from the definition $\tilde{y}_j:=\hat{\theta}_jy_j+(1-\hat{\theta}_j)(\hat{B}^{(0)}_{j+1}+\delta I)s_j-\gamma s_j$ with $\hat{\theta}_j$ given in (\ref{damp3}), it follows from \textit{Lemma 1} that inequalities $s^T_j\tilde{y}_j\geq 0.2s^T_j(B^{(0)}_{j+1}+\delta I)s_j>0$ hold. This yields:
\begin{equation}\label{proofeq3}
\begin{aligned}
\frac{\tilde{y}^T_{j}\tilde{y}_{j}}{s^T_{j}\tilde{y}_{j}}&\leq\frac{\norm{\hat{\theta}_jy_j+(1-\hat{\theta}_j)(\hat{B}^{(0)}_{j+1}+\delta I)s_j-\gamma s_j}^2}{0.2s^T_j(B^{(0)}_{j+1}+\delta I)s_j}\\
&=\frac{1}{0.2s^T_j(B^{(0)}_{j+1}+\delta I)s_j} \{\hat{\theta}^2_jy^T_jy_j+(1-\hat{\theta}_j)^2\\
&\quad\quad \cdot s^T_j(\hat{B}^{(0)}_{j+1}+\delta I)^2s_j+2\hat{\theta}_j(1-\hat{\theta}_j)\\
&\quad\quad\cdot y^T_j(\hat{B}^{(0)}_{j+1}+\delta I)s_j+\gamma^2s^T_js_j\\
&\quad\quad-2\gamma s^T_j[\hat{\theta}_jy_j+(1-\hat{\theta}_j)(\hat{B}^{(0)}_{j+1}+\delta I)s_j]\}  .  
\end{aligned}	
\end{equation}
From the definition $y_j=\frac{1}{m_j}\sum_{l=1}^{m_j}\nabla{F(\bar{x}_{j+1},\xi_{j,l})}-\nabla{F(\bar{x}_{j},\xi_{j,l})}$, and using the first-order Taylor approximation at $\bar{x}_j$, we have   $y_j=\frac{1}{m_j}\sum_{l=1}^{m_j}\nabla^2{F(\bar{x}_{j}+\vartheta s_j,\xi_{j,l})}s_j$, where $0<\vartheta<1$. Thus,  $y^T_jy_j=\frac{1}{m^2_j}s^T_j\{\sum_{l=1}^{m_j}\sum_{r=1}^{m_j}\nabla^2{F(\bar{x}_{j}+\vartheta s_j,\xi_{j,r})}\nabla^2{F(\bar{x}_{j}+\vartheta s_j,\xi_{j,l})}\}s_j$. With Assumption 1 that $\norm{\nabla^2F(x,\xi)}\leq \rho$, which implies $-\rho I\prec\nabla^2F(x,\xi)\prec \rho I$. We further have $y^T_jy_j\leq \rho^2s^T_js_j$. Next, we consider the product $y^T_js_j$. Since $y^T_js_j=\frac{1}{m_j}\sum_{l=1}^{m_j}s^T_j\nabla^2{F(\bar{x}_{j}+\vartheta s_j,\xi_{j,l})}s_j$, it follows that $-\rho s^T_js_j \leq y^T_js_j \leq \rho s^T_js_j$.  Substituting the above inequality into (\ref{proofeq3}), with $\hat{B}^{(0)}_{j+1}=\tau_{j+1}I$, we get:
\begin{equation}\label{proofeq4}
\begin{aligned}
\frac{\tilde{y}^T_{j}\tilde{y}_{j}}{s^T_{j}\tilde{y}_{j}}&\leq \frac{1}{0.2(\tau_{j+1}+\delta)}\{\hat{\theta}^2_j\rho^2+(1-\hat{\theta}_j)^2(\tau_{j+1}+\delta)^2\\
&\quad\quad+2\rho\hat{\theta}_j(1-\hat{\theta}_j)(\tau_{j+1}+\delta)+\gamma^2\\
&\quad\quad+2\gamma\hat{\theta}_j\rho-2(1-\hat{\theta}_j)(\tau_{j+1}+\delta)\gamma\}\\
&=\frac{5(\hat{\theta}_j\rho+\gamma)^2}{\tau_{j+1}+\delta}+5(1-\hat{\theta}_j)^2(\tau_{j+1}+\delta)\\
&\quad\quad +10\hat{\theta}_j(1-\hat{\theta}_j)\rho-10\gamma(1-\hat{\theta}_j).
\end{aligned}
\end{equation}   
By using $\tau_{j+1}=\text{max}\left\{\frac{y^T_jy_j}{s^T_jy_j}+\gamma,\beta \right\}$, we have $\beta+\delta\leq\tau_{j+1}+\delta\leq\beta+\rho+\gamma+\delta$. Furthermore, $10\hat{\theta}_j(1-\hat{\theta}_j)\rho\leq 5\rho(1-\hat{\theta}^2_j)$ holds true as $0<\hat{\theta}_j\leq 1$. By using the property of the function $\varphi(x)=ax+\frac{b}{x},\;a>0,b>0$, we obtain the following result:
\begin{equation}\label{proofeq5}
\begin{aligned}
\frac{\tilde{y}^T_{j}\tilde{y}_{j}}{s^T_{j}\tilde{y}_{j}}&\leq Q+10\hat{\theta}_j(1-\hat{\theta}_j)\rho-10\gamma(1-\hat{\theta}_j)\\
&\leq Q+5\rho(1-\hat{\theta}^2_j)-10\gamma(1-\hat{\theta}_j)\leq Q+5\rho,
\end{aligned}	
\end{equation}  
where $Q$ is defined in (\ref{Q}). Therefore, by substituting the results in (\ref{proofeq5}) into (\ref{proofeq2}), one gets $\norm{\hat{B}^{(i+1)}_{t}}\leq\norm{\hat{B}^{(i)}_{t}}+Q+5\rho+\gamma$. By induction, we then obtain the desired result:
\begin{equation}\label{proofeq6}
\norm{\hat{B}^{(\tilde{M})}_{t}}\leq \norm{\hat{B}^{(0)}_{t}}+\tilde{M}(Q+5\rho+\gamma)\leq\beta+\rho+\gamma+\tilde{M}(Q+5\rho+\gamma).
\end{equation}

Thus, we have proved the upper bound on the norm of the matrix $\hat{B}^{(\tilde{M})}_{t}$, the next lemma gives for a more accurate lower bound rather than just $\hat{B}^{(\tilde{M})}_{t}\succeq \gamma I$.

\textit{Lemma 3.}  Given the initial positive definite matrix $\hat{B}^{(0)}_t$ defined by (\ref{B}), and suppose $\hat{B}^{(i+1)}_t$ is updated via L-BFGS step of Algorithm 1, then with Assumption 1, all eigenvalues  of $\hat{B}^{(\tilde{M})}_t$ satisfies
\begin{equation}\label{B_belowBound}
\lambda(\hat{B}^{(\tilde{M})}_t)\geq Q_L,
\end{equation} 
where $	Q_L=\text{max}\left\{\tilde{Q}^{-1},\gamma^{-1}\right\}$ and 
\begin{equation}
\begin{aligned}
\tilde{Q}=&\frac{w^{2\tilde{M}}-1}{Q+5\rho+2\sqrt{0.2(Q+5\rho)(\beta+\delta)}}+\beta^{-1}w^{2\tilde{M}},
\end{aligned}
\end{equation}
with $w:=\sqrt{\frac{Q+5\rho}{0.2(\beta+\delta)}} +1$.

\textit{Proof.} From (\ref{res3}), we have:
\begin{equation}\label{proofeq7}
\hat{B}^{(i+1)}_{t}\succeq\hat{B}^{(i)}_{t}+\frac{\tilde{y}_{j}\tilde{y}^T_{j}}{s^T_{j}\tilde{y}_{j}}-\frac{\hat{B}^{(i)}_{t}s_{j}s^T_{j}\hat{B}^{(i)}_{t}}{s^T_{j}\hat{B}^{(i)}_{t}s_{j}}.
\end{equation}
Since both sides of the inequality (\ref{proofeq7}) are positive definite matrices, taking matrix inversion and using the Sherman–
Morrison–Woodbury formula yields:
\begin{equation}\label{proofeq8}
\begin{aligned}
\hat{H}^{(i+1)}_t&\preceq \left(I-\frac{s_j\tilde{y}^T_j}{s^T_j\tilde{y}_j}\right)\hat{H}^{(i)}_t\left(I-\frac{\tilde{y}_js^T_j}{s^T_j\tilde{y}_j}\right)+\frac{s_js^T_j}{s^T_j\tilde{y}_j}\\
&=\hat{H}^{(i)}_t-\frac{1}{s^T_j\tilde{y}_j}(s_j\tilde{y}^T_j\hat{H}^{(i)}_t+\hat{H}^{(i)}_t\tilde{y}_js^T_j)+\frac{\tilde{y}^T_j\hat{H}^{(i)}_t\tilde{y}_j}{(s^T_j\tilde{y}_j)^2}\\
&\quad\quad\quad\quad\cdot s_js^T_j+\frac{s_js^T_j}{s^T_j\tilde{y}_j},
\end{aligned}	
\end{equation}
where $\hat{H}^{(i)}_t$ is the  inverse matrix of $\hat{B}^{(i)}_t$, i.e., $\hat{H}^{(i)}_t:=\hat{B}^{(i)^{-1}}_t$. By taking  the matrix norm on both sides of (\ref{proofeq8}) and using the triangle inequality, we get:
\begin{equation}\label{proofeq9}
\begin{aligned}
\norm{\hat{H}^{(i+1)}_t}\leq&\norm{\hat{H}^{(i)}_t}+\frac{2\norm{\hat{H}^{(i)}_t}\cdot\norm{s_j}\cdot\norm{\tilde{y}_j}     }{s^T_j\tilde{y}_j}+\frac{\tilde{y}^T_j\tilde{y}_j  }{s^T_j\tilde{y}_j}\cdot\frac{s^T_js_j}{s^T_j\tilde{y}_j}\\
&\cdot\norm{\hat{H}^{(i)}_t}+\frac{s^T_js_j}{s^T_j\tilde{y}_j}.
\end{aligned}
\end{equation}
Recall from  the proof of \textit{Lemma 2} that $\frac{\tilde{y}^T_j\tilde{y}_j  }{s^T_j\tilde{y}_j}\leq Q+5\rho$. Moreover, according to \textit{Lemma 1}, we have  $\frac{s^T_js_j}{s^T_j\tilde{y}_j}\leq \frac{s^T_js_j}{0.2s^T_j(\hat{B}^{(0)}_{j+1}+\delta I) s_j}=\frac{1}{0.2(\tau_{j+1}+\delta)}$ and hence
\begin{equation}
\frac{\norm{s_j}\cdot\norm{\tilde{y}_j}     }{s^T_j\tilde{y}_j}=\left(\frac{s^T_js_j}{s^T_j\tilde{y}_j}\cdot\frac{\tilde{y}^T_j\tilde{y}_j}{s^T_j\tilde{y}_j}\right)^{1/2}\leq \sqrt{\frac{Q+5\rho}{0.2(\tau_{j+1}+\delta)}}.
\end{equation} 
Substituting the above results into (\ref{proofeq9}) and noting the fact $\tau_{j+1}\geq \beta$, (\ref{proofeq9}) can be simplified to
\begin{equation}\label{proofeq10}
\norm{\hat{H}^{(i+1)}_t}\leq w^2\norm{\hat{H}^{(i)}_t}+\frac{1}{0.2(\beta+\delta)}. 
\end{equation}
By induction with $\hat{H}^{(0)}_t\preceq\beta^{-1} I$, we obtain the desired result.

Based on the above uniformly upper bound and lower bound on the resultant L-BFGS matrix, we now derive the convergence result of our proposed algorithm. Moreover, the following assumption is required.

\textit{Assumption 2}. For any iteration, the variance of the gradient conditioned on current iterate is bounded above:
\begin{equation}\label{converge1}
\mathbb{E}(\norm{\nabla F(x_k,\xi_k)-\nabla f(x_k)}^2|x_k)\leq \sigma^2.
\end{equation}
Moreover, the norm square of the  gradient is expected to be bounded above by a positive constant $D$ \cite{oLBFGS,MokhtariAryan2014RRSB,RHBYRD}:
\begin{equation}\label{converge1_1}
\mathbb{E}[\norm{\nabla F(x_k,\xi_k)}^2|x_k]\leq D.
\end{equation}

With Assumption 2, we introduce the following lemma:

\textit{Lemma 4}\cite{SQN,MokhtariAryan2014RRSB}. Suppose Assumption 2 holds, and the sequence $\{x_k\}$ for $k=1,\dots,$ is generated with the initial value $x_0$ and using a specific constant batch size $m_k=m$. Then there exists a positive constant $M_f$ such that $\mathbb{E}[f(x_k)]\leq M_f$. Moreover, the sequence almost surely converges to a stationary point, i.e., $\lim\limits_{\text{$k\rightarrow\infty$}}\norm{\nabla f(x_k)}=0,\;\text{with probability} \; 1.$

We are now ready to proceed to show the convergence of our proposed algorithm under the given assumptions, which is summarized in the following theorem. Without loss of generality, the interval length is assumed to be unity.

\textit{Theorem 1}. Suppose the iterations of the Sd-REG-LBFGS algorithm satisfies Assumption 2, and the sequence $\{x_k\}$ for $k=1,\dots,N-1$ is generated with initial value $x_0$. Given the constant batch size $m_k=m$ and in particular the following step size:
\begin{equation}\label{converge2}
\eta_k=\frac{\eta_0Q^{-1}_U}{k^\upsilon+(L_f/2)\eta_0Q^{-2}_L},
\end{equation} 
with $0.5<\upsilon<1$, the following inequality holds:
\begin{equation}\label{converge3}
\begin{aligned}
\frac{1}{N}\sum_{k=0}^{N-1}\mathbb{E}&\left(\norm{\nabla f(x_k)}^2 \right)\leq \frac{ [(N-1)^\upsilon+(L_f/2)\eta_0Q^{-2}_L]^2 }{\eta_0Q^{-2}_U(N-1)^\upsilon N}\\
&\cdot(M_f-f^{l})+\frac{L_fQ^{-2}_L\sigma^2\eta_0[(N-1)^{1-\upsilon}-1]}{2m(1-\upsilon)N},
\end{aligned}
\end{equation}
where  $f^l:=\text{min}\{f(x_0),\dots,f(x_{N-1})\}$ and $N$ is the iteration number. Furthermore, given a constant $0<\epsilon<1$, the iteration number $N$ needed to ensure $\frac{1}{N}\sum_{k=0}^{N-1}\mathbb{E}\left(\norm{\nabla f(x_k)}^2 \right)\leq \epsilon$ is at most $O(\epsilon^{-\frac{1}{1-\upsilon}})$.

\textit{Proof}. Recall that the gradient of $f(\cdot)$ is Lipschitz continuous with constant $L_f$, therefore, using second-order Taylor expansion at iteration $k$ leads to:
\begin{equation}\label{converge4}
\begin{aligned}
&f(x_{k+1})\leq f(x_k)+\nabla f(x_k)^T(x_{k+1}-x_k)+\frac{L_f}{2}\norm{x_{k+1}-x_k}^2\\
&\quad\quad=f(x_k)+\nabla f(x_k)^T(-\eta_k\hat{B}^{-1}_k\bar{g}_k)+\frac{L_f}{2}\eta^2_k\norm{\hat{B}^{-1}_k\bar{g}_k}^2 \\
&\quad\quad\leq f(x_k)-\eta_k\nabla f(x_k)^T\hat{B}^{-1}_k\bar{g}_k+\frac{L_f}{2}\eta^2_k\norm{\hat{B}^{-1}_k}^2\cdot \norm{\bar{g}_k}^2,
\end{aligned}
\end{equation} 
where for notational convenience, we denote $\bar{g}_k=\bar{g}_k(x_k,\xi_k)$. From Lemma 2 and Lemma 3, we have $Q^{-1}_UI\preceq\hat{B}^{-1}_k\preceq Q^{-1}_LI$. Substituting it into (\ref{converge4}) results in:
\begin{equation}\label{converge4_1}
f(x_{k+1})\leq f(x_k)-\eta_kQ^{-1}_U\nabla f(x_k)^T\bar{g}_k+\frac{L_f}{2}\eta^2_kQ^{-2}_L\norm{\bar{g}_k}^2.
\end{equation}
To evaluate the expectation of (\ref{converge4_1}), we shall first take the expectation conditioned on $x_k$ on both sides and then the expectation with respect to $x_k$. We shall make use of the fact that  $\mathbb{E}_B[\mathbb{E}_A(A|B)]=\mathbb{E}(A)$ for random variables $A$ and $B$. Consequently, with Assumption 2, we get:
\begin{equation}\label{converge5}
\begin{aligned}
\mathbb{E}[f(x_{k+1})]
&\leq \mathbb{E}[f(x_k)]-\eta_kQ^{-1}_U\mathbb{E}[\nabla f(x_k)^T\mathbb{E}(\bar{g}_k|x_k)]\\
&\quad\quad+\frac{L_f}{2}\eta^2_kQ^{-2}_L\mathbb{E}[\mathbb{E}(\norm{\bar{g}_k}^2|x_k)].
\end{aligned}
\end{equation}
Furthermore, we have:
\begin{equation}\label{converge6}
\begin{aligned}
\mathbb{E}(\norm{\bar{g}_k-\nabla f(x_k)}^2|x_k)&=\mathbb{E}(\norm{\bar{g}_k}^2|x_k)-\norm{\nabla f(x_k)}^2,
\end{aligned}	
\end{equation}
and it further yields $\mathbb{E}(\norm{\bar{g}_k}^2|x_k)=\sigma^2/m+\norm{\nabla f(x_k)}^2$.
Substituting the result into (\ref{converge5}), we have:
\begin{equation}\label{converge8}
\begin{aligned}
\mathbb{E}[f(x_{k+1})]&\leq\mathbb{E}[f(x_k)]-\left(\eta_kQ^{-1}_U-\frac{L_f}{2}\eta^2_kQ^{-2}_L\right)\\
&\cdot\mathbb{E}(\norm{\nabla f(x_k)}^2)+\frac{L_f\eta^2_kQ^{-2}_L\sigma^2}{2m}.
\end{aligned}
\end{equation}
By summing (\ref{converge8}) for $k=0,\dots,N-1$, the following result is obtained: 
\begin{equation}\label{converge9}
\begin{aligned}
\sum_{k=0}^{N-1}\mathbb{E}&(\norm{\nabla f(x_k)}^2)\leq\sum_{k=0}^{N-1}\frac{\mathbb{E}[f(x_{k})]-\mathbb{E}[f(x_{k+1})]}{\eta_kQ^{-1}_U-(L_f/2)\eta^2_kQ^{-2}_L}\\
&\quad+\sum_{k=0}^{N-1}\frac{L_f\eta_kQ^{-2}_L\sigma^2}{2m[Q^{-1}_U-(L_f/2)\eta_kQ^{-2}_L]}.
\end{aligned}
\end{equation}
Furthermore, from (\ref{converge2}), we have $\frac{\eta_k}{Q^{-1}_U-(L_f/2)\eta_kQ^{-2}_L}=\eta_0k^{-\upsilon}$. 
Substituting it into (\ref{converge9}), we obtain the simplified inequality:
\begin{equation}\label{converge11}
\begin{aligned}
\sum_{k=0}^{N-1}\mathbb{E}&(\norm{\nabla f(x_k)}^2)\leq\sum_{k=0}^{N-1}\frac{\eta_0k^{-\upsilon}}{\eta^2_k}(\mathbb{E}[f(x_{k})]\\
&-\mathbb{E}[f(x_{k+1})])+\frac{L_fQ^{-2}_L\sigma^2\eta_0}{2m}\sum_{k=0}^{N-1}k^{-\upsilon}.
\end{aligned}
\end{equation}
By utilizing the result in Lemma 3 that  $\mathbb{E}[f(x_k)]\leq M_f$, we further have:
\begin{equation}\label{converge12}
\begin{aligned}
&\sum_{k=0}^{N-1}\mathbb{E}(\norm{\nabla f(x_k)}^2)\leq \sum_{k=1}^{N-1}\left(\frac{\eta_0k^{-\upsilon}}{\eta^2_k}-\frac{\eta_0(k-1)^{-\upsilon}}{\eta^2_{k-1}}\right)\mathbb{E}[f(x_{k})]\\
&\quad-\frac{\eta_0(N-1)^{-\upsilon}}{\eta^2_{N-1}}\mathbb{E}[f(x_{N})]+\frac{L_fQ^{-2}_L\sigma^2\eta_0}{2m}\sum_{k=0}^{N-1}k^{-\upsilon}\\
&\quad\leq M_f\sum_{k=1}^{N-1}\left(\frac{\eta_0k^{-\upsilon}}{\eta^2_k}-\frac{\eta_0(k-1)^{-\upsilon}}{\eta^2_{k-1}}\right)\\
&\quad-\frac{\eta_0(N-1)^{-\upsilon}}{\eta^2_{N-1}}f^l+\frac{L_fQ^{-2}_L\sigma^2\eta_0}{2m}\sum_{k=0}^{N-1}k^{-\upsilon}\\
&\quad= \frac{\eta_0(M_f-f^l)(N-1)^{-\upsilon}}{\eta^2_{N-1}}+\frac{L_fQ^{-2}_L\sigma^2\eta_0}{2m}\sum_{k=0}^{N-1}k^{-\upsilon}\\
&\quad=\frac{ [(N-1)^\upsilon+(L_f/2)\eta_0Q^{-2}_L]^2(M_f-f^{l}) }{\eta_0Q^{-2}_U(N-1)^\upsilon }\\
&\quad\quad\quad\quad\quad+\frac{L_fQ^{-2}_L\sigma^2\eta_0}{2m}\sum_{k=0}^{N-1}k^{-\upsilon}.
\end{aligned}	
\end{equation}
By applying following inequality:
\begin{equation}\label{converge12_1}
k^{-\upsilon}\leq \frac{k^{1-\upsilon}-(k-1)^{1-\upsilon}}{1-\upsilon},\;\text{for} \;k\geq 1,
\end{equation}
to (\ref{converge12}), we obtain the desired result in (\ref{converge3}). For a given constant $\epsilon$ satisfying $0<\epsilon<1$, the iteration number needed to guarantee $\frac{1}{N}\sum_{k=0}^{N-1}\mathbb{E}(\norm{\nabla f(x_k)}^2)<\epsilon$ satisfies:
\begin{equation}\label{converge13}
\begin{aligned}
&\frac{ [(N-1)^\upsilon+(L_f/2)\eta_0Q^{-2}_L]^2(M_f-f^{l}) }{\eta_0Q^{-2}_U(N-1)^\upsilon N}+\\
&\quad\quad\quad\frac{L_fQ^{-2}_L\sigma^2\eta_0[(N-1)^{1-\upsilon}-1]}{2m(1-\upsilon)N}<\epsilon.
\end{aligned}
\end{equation}
Therefore, for $0.5<\upsilon<1$, the iteration number is at most $O(\epsilon^{-\frac{1}{1-\upsilon}})$ to reach $\frac{1}{N}\sum_{k=0}^{N-1}\mathbb{E}(\norm{\nabla f(x_k)}^2)<\epsilon$.

\section{Empirical Study}
We have studied the theoretical properties and the convergence of the proposed quasi-Newton method in the previous section. In this section, we will apply the proposed method to solve several  optimization problems in machine learning. Specifically, two machine learning problems will be studied, namely logistic regression and Bayesian logistic regression for binary classification. To carry out the optimization, the gradient required by the Sd-REG-LBFGS method is obtained analytically. In a general fashion, we mainly focus on nonconjugate exponential models under stochastic regime, in which Bayesian logistic regression is a particular example. 
\subsection{Logistic Regression}
We first consider the logistic regression problem. The objective function is given as follows \cite{Bishop}:
\begin{equation}\label{logreg}
f(\theta)=-\frac{1}{N}\sum_{n=1}^{N}z_n\text{log }\sigma(\theta^Tx_n)+(1-z_n)\text{log }\sigma(-\theta^Tx_n),
\end{equation}  
where $\sigma(\cdot)$ is the sigmoid function given by $\sigma(x)=1/(1+\text{exp}(-x))$, $x_n$ is the feature vector and $z_n$ is its label. 
\subsection{Sd-REG-LBFGS for VBI}
Variational Bayesian inference (VBI) is an efficient method for approximating the a posteriori probability distributions for making inference. The main ingredient is to convert inference problems into optimization problems with the KL-divergence as the objective function. Another popular scheme for making inference is Markov chain Monte Carlo (MCMC) sampling method. It can be easily parallelized for multiple processors to reduce the computational cost for high dimension problems. In this section, we illustrate the application of the proposed Sd-REG-LBFGS to the delta VBI scheme for nonconjugate models proposed in \cite{NONCON}. The resultant algorithm is denoted by SDVBI. In addition, interested readers can refer to \cite{ARCG,HoffmanSVI,VBISS} for applications of optimization methods in VBI.

Suppose $x_{1:N}$ are observations,  $z_{1:N}$ are local hidden variables and $\theta$ is global hidden variable. Furthermore, $\theta$ is the nonconjugate variable and $z_{1:N}$ are conjugate variables. Consider the nonconjugate model in \cite{NONCON} as follows: 
\begin{equation}\label{vb2}
p(x,z,\theta)=p(\theta)\cdot\prod_{n=1}^{N}p(x_n|z_n)p(z_n|\theta),
\end{equation}
where  $p(z_n|\theta)=h(z_n)\text{exp}\{\eta_{g_n}(\theta)^Tt(z_n)-a(\eta_{g_n}(\theta))\}$ and  $p(x_n|z_n)=h(x_n)\text{exp}\{ t(z_n)^T[t(x_n)^T,1]^T \}$.
The goal of variational inference is to approximate the posterior distribution by  finding a member of a specific family $\mathcal{Q}$ to minimize its KL-divergence to the true a posteriori distribution:   
\begin{equation}\label{vb3}
q^*(z,\theta)=\operatorname{argmin}_{q\in{\mathcal{Q}}}\: \text{KL}(q(z,\theta)|| p(z,\theta|x)).
\end{equation}

For the MFVI framework, the statistical independence between hidden variables with a fully factorized variational distribution family are assumed, i.e.,
\begin{equation}\label{vb4}
q(z,\theta)=q(\theta|\lambda)\cdot\prod_{n=1}^{N}q(z_n|\varphi_n),
\end{equation}
where $q(z_n|\varphi_n)=h(z_n)\text{exp}\{\eta_l(\varphi_n)^Tt(z_n)-a(\eta_l(\varphi_n))\},$ and  Gaussian distribution has been adopted to approximate its variational distribution, i.e., $q(\theta|\lambda)=\mathcal{N}(\mu,S)$, with $\lambda$ being the parameter pair $(\mu,S)$. The following can be obtained by substituting the results into (\ref{vb3}):
\begin{equation}\label{vb5}
\begin{aligned}
&\text{KL}(q||p)=\mathbb{E}_q[\text{log}\:q(z_{1:N},\theta)]-\mathbb{E}_q[\text{log}\:p(z_{1:N},\theta|x_{1:N})]\\
&\quad\quad=\mathbb{E}_q[\text{log}\:q(z_{1:N},\theta)]-\mathbb{E}_q[\text{log}\:p(x_{1:N},z_{1:N},\theta)]+const.\\
&\quad\quad:=\mathcal{L}(q).
\end{aligned}
\end{equation}
First, for nonconjugate variable $\theta$, the objective function of delta VBI has been derived based on second-order Taylor approximation of the variational objective function in \cite{NONCON} as follows:
\begin{equation}\label{vb12}
\begin{aligned}
\mathcal{L}(\lambda)&=\mathbb{E}_{q(\theta|\lambda)}[\text{log}\;q(\theta|\lambda)]-\sum_{n=1}^{N}\mathbb{E}_{q(\theta,z_n)}[\text{log}\;p(z_n|\theta)]\\
&\approx d(\mu)+\frac{1}{2}(\text{Tr}\{\nabla^2d(\mu)S\}-\text{log det}S)+const.,
\end{aligned}
\end{equation}
where $d(\theta):= -\eta_g(\theta)^T\cdot\sum_{n=1}^{N}\nabla_{\eta_l}a(\eta_l(\varphi_n))+Na(\eta_g(\theta))-\text{log}\;p(\theta)$, the optimization problem becomes:
\begin{equation}\label{sla5}
\lambda^*=\operatorname{argmin}_{\lambda\in{\mathbb{R}^d}}\;\{\mathcal{L}(\lambda)=-\frac{1}{2}\text{log det}S+\mathbb{E}_{q(\theta|\lambda)}d(\theta)\}.
\end{equation} 
We notice that (\ref{vb12}) contains large summation term, which makes the gradient evaluation computationally rather expensive. Next, we randomly sample a subset $\mathcal{S}$ from $\{1,\dots,N\}$ to form an unbiased stochastic gradient, which is denoted as $\nabla_{\lambda}\mathcal{L}(\lambda;\mathcal{S})$. We omit the reduplicative and tedious derivation, as the full gradient can be found in Appendix C of \cite{NONCON}.

For the conjugate variable $z_n$ updating, the variational objective function $\mathcal{L}(\varphi_n)$ from the KL-divergence in (\ref{vb5}) in \cite{NONCON} is as follows:
\begin{equation}\label{vb19}
\begin{aligned}
\mathcal{L}(\varphi_n)&=\mathbb{E}_{q(z_n)}[\text{log}\;q(z_n)]-\mathbb{E}_{q(z_n)}[\text{log}\;p(x_n|z_n)]\\
&\quad\quad-\mathbb{E}_{q(z_n,\theta)}[\text{log}\;p(z_n|\theta)]+const\\
&=\{\eta_l(\varphi_n)^T-[t(x_n)^T,1]-\mathbb{E}_{q(\theta)}[\eta_g(\theta)^T] \}\\
&\quad\quad\cdot\nabla_{\eta_l}a(\eta_l(\varphi_n))-a(\eta_l(\varphi_n))+const,
\end{aligned}
\end{equation}
where the last equality in (\ref{vb19}) follows from the basic property of the exponential family. To derive the update for $\varphi_n$, we take the gradient of $\mathcal{L}(\varphi_n)$:
\begin{equation}\label{vb20}
\begin{aligned}
\nabla_{\varphi_n}\mathcal{L}(\varphi_n)&=\mathcal{D}_{\varphi_n}\eta_l(\varphi_n)^T\cdot\nabla^2_{\eta_l}a(\eta_l(\varphi_n))\{\eta_l(\varphi_n)\\
&\quad\quad-[t(x_n)^T,1]^T-\mathbb{E}_{q(\theta)}[\eta_g(\theta)] \},
\end{aligned}
\end{equation} 
where $\mathcal{D}_{\varphi_n}\eta_l(\varphi_n)$ is the Jacobian matrix of $\eta_l(\cdot)$ with respect to $\varphi_n$. Therefore, by using the gradient for optimization or by simply setting the gradient to zero, i.e., $\nabla_{\varphi_n}\mathcal{L}(\varphi_n)=0$, we  obtain the conjugate variable update. With the above stochastic gradients derived, we have shown the application of the proposed method.

In particular, with the following settings \cite{NONCON}:
\begin{equation}
\begin{aligned}
&h(z_n)=1,\;t(z_n)=[z_n,1-z_n]^T, a(\eta_g(\theta))=0,\\
&\eta_{g_n}(\theta)=[\text{log }\sigma(\theta^Tx_n),\text{log }\sigma(-\theta^Tx_n)]^T,\; n=1,\dots,N,
\end{aligned}
\end{equation}
one recovers Bayesian logistic regression. Here, it should be noted that VBI is only considered for the nonconjugate variable $\theta$. However, for the settings of correlated topic model, VBI is considered for both $\theta$ and $z_n$. As the applications of Sd-REG-LBFGS are similar, we shall consider Bayesian logistic regression for numerical experiments for simplicity.

\section{Numerical Results}

In this section, the numerical experiments are performed on our proposed Sd-REG-LBFGS algorithm. Two applications are considered in machine learning, which are logistic regression and SDVBI for Bayesian logistic regression. Moreover, in this paper, we only consider binary classification problems. We also employed a synthetic dataset and several real datasets \cite{scene,Elec,wifi1,wifi2,IONO,BNA} for the performance evaluation. In particular for the parameter studies, we use a synthetic dataset and a real \textit{scene} dataset \cite{scene} (available at http://mulan.sourceforge.net /datasets-mlc.html), which  can be categorized as the following 4 scenarios:

S1. Solving logistic regression (LR) using synthetic dataset, which is presented in Section V-A;

S2. Solving Bayesian logistic regression (BLR) using synthetic dataset, which is presented in Section V-A;

S3.  Solving LR using \textit{scene} dataset in \cite{scene}. Due to page limitation, the results are presented in Section II of the supplementary material;

S4. Solving BLR  using \textit{scene} dataset in \cite{scene}. Due to page limitation, the results are presented in Section II of the supplementary material.\\

The following algorithms are considered for evaluation:

(A) Proposed Sd-REG-LBFGS:  The proposed stochastic damped regularized L-BFGS as described in Algorithm 1;

(B) SdLBFGS: Stochastic damped regularized L-BFGS without regularization in \cite{SQN};

(C) SGD: Stochastic gradient descent is adopted;

(D) SAA: Stochastic approximation averaging in [39] is applied;

(E) RSA:  Robust stochastic approximation in \cite{RobustSA} is employed.

(F) Adam: Adam \cite{adam} is employed. 

Here, we summarize again some key parameters and variables that are involved in the numerical experiments.  

1. $d$: the dimension of the optimization variable, e.g., $\theta\in\mathbb{R}^d$. 

2. $N$: the number of training points in the dataset. 

3. $m$: the batch size used for stochastic approximation of the gradient. In the numerical experiment, we use constant batch size at each iteration, e.g., $m=|\mathcal{S}|$  for $\nabla_{\theta}\mathcal{L}(\theta;\mathcal{S})$. 

4. $M$: the memory size used for the Sd-REG-LBFGS algorithm to store the correction pairs $(s_t,y_t)$  calculated by (20) and (21).

5. $L$: the interval length. Every $L$  iterations, we perform averaging on the   iterate points, which is used to calculate the correction pairs by (20) and (21). 

6. $\gamma$: the regularized parameter for BFGS update given in (24), which prevent the L-BFGS matrix from being close to singularity. 

7. $\eta_k$: the step size for SGD, SdLBFGS and Sd-REG-LBFGS optimization schemes. In the numerical experiment, we adopt the diminishing step size $\eta_k=r/k$  with a positive constant $r$  at each iteration.

In general, for the regression problems, one needs to include a constant bias term. This can be implemented by concatenating a unity element at the beginning or the end of each input vector, i.e., if the unity is put at the beginning, $\theta_0+\theta^Tx_n=[\theta_0,\theta^T][1,x^T_n]$. For notational convenience, we omit the bias term here. The performance of various algorithms will be evaluated in terms of the norm of the gradient (NOG) and the classification accuracy (ACC). The NOG for LR is defined as follows: 
\begin{equation}\label{NOG1}
	\text{NOG}=\norm{\frac{1}{N}\sum_{n=1}^{N}[z_n-\sigma(\theta^Tx_n)]x_n}.
\end{equation}
Moreover, the exact gradient of the objective function in BLR can be calculated as follows \cite{NONCON}:
\begin{equation}\label{NOG2}
	\begin{aligned}
	\nabla_{\theta}\mathcal{L}(\theta)&=\frac{1}{N}\sum_{n=1}^{N}\{[z_n-\sigma(\theta^Tx_n)]x_n+\frac{1}{2}\sigma(\theta^Tx_n)\\
	&\cdot\sigma(-\theta^Tx_n)[1-2\sigma(\theta^Tx_n)]x_nx^T_nSx_n\}+S^{-1}_0.
	\end{aligned}
\end{equation}
Hence, the NOG for BLR is defined as  $\text{NOG}=\norm{\nabla_{\theta}\mathcal{L}}$. 

Lower NOG indicates the better convergence of an algorithm to a stationary point. The classification accuracy is given as
\begin{equation}
	ACC=\frac{TP+TN}{TP+FN+FP+FN}
\end{equation}
where TP, TN, FP, and FN denote true positives, true negatives, false positives and false negatives, respectively.  The decision rules for class prediction are given as 
\begin{equation}
	\text{if } \sigma(\hat{\theta}^Tx_n)\geq 0.5, \text{ then } z_n=1, \text{ else } z_n=0,	
\end{equation}
for logistic regression and Bayesian logistic regression respectively, where $\hat{\theta}$ is the estimated value of $\theta$. We have also conducted a sensitivity analysis to study the effects of different batch sizes, memory sizes and regularization parameters with our proposed method. Due to page limitation, the details are omitted here. Interested readers are referred to Section II of the supplementary material for details.  

\subsection{Performance comparisons with different real datasets}
In this subsection, we first study the effectiveness of our proposed method using the same settings of each algorithms for their common parameters with different real datasets. Specifically, for all schemes, batch size is set to a relatively small value to show that our proposed method is particularly effective. Moreover, the real datasets are described as follows:

1. \textit{Banknote Authentication Dataset} (BNA) \cite{BNA} (available at UCI Machine Learning Repository): we use 1,370 samples, which has 4 variables. Considering 5-fold cross validation, there are 1,096 data points for training and 274 samples for testing. 

2. \textit{Wireless Indoor Localization Dataset} (WINL) \cite{wifi1,wifi2} (available at UCI Machine Learning Repository): 2,000 samples with 7 variables are used for the performance evaluation. For 5-fold cross-validation, there are 1,600 data points for training and 400 samples for testing. 

3. \textit{Ionosphere Dataset} (IONO) \cite{IONO} (available at UCI Machine Learning Repository): we use 350 samples with 33 variables for performance evaluation. There are 280 data points for training and 70 samples for testing according to 5-fold cross validation.

4. \textit{Electrical Grid Stability Simulated Dataset} (ELEG) \cite{Elec} (available at UCI Machine Learning Repository):  10,000 samples of the dataset with 14 variables is used for performance evaluation, of which 8,000 are for training and 2,000 are for testing according to 5-fold cross validation.

Moreover, the batch size and step size for each algorithm are set to $m=20$ and $\eta_k=7/k$, respectively. For our proposed method and SdLBFGS, the memory is set to the same value $M=10$. We set the regularization parameters for our proposed method to $\gamma=10^{-4}$ and $\delta=1.25+0.01$, respectively. For each NOG and ACC value, it is computed via the average of 5-fold cross validation and 50 Monte Carlo runs. The results are shown in Table \ref{table_datasets}. It can be seen that our proposed method performs the best obviously in terms of NOG performance. For ACC performance evaluation, our proposed method is generally better than other methods, except that the proposed method is slightly worse than SdLBFGS for WINL and ELEG. This is due to the bias that our method has introduced. However, Sd-REG-LBFGS is more robust as SdLBFGS has resulted in ill-conditioning problems during the experiments. 

Next, we will consider the synthetic dataset and the real dataset \textit{scene} \cite{scene} to extensively study  the effects of different parameter settings.

\newcolumntype{P}[1]{>{\centering\arraybackslash}p{#1}}
\newcolumntype{M}[1]{>{\centering\arraybackslash}m{#1}}
\begin{table}[!t]
	\begin{center}
		\renewcommand{\arraystretch}{1.3}
		\caption{The NOG and ACC performances of various algorithms averaged over 50 Monte Carlo simulations and 5-fold cross validation with different datasets for logistic regression and Bayesian logistic regression.}
		\label{table_datasets}
		\centering
		
		\begin{tabular}{|M{0.9cm}|M{1.2cm}|M{1cm}|M{1cm}|M{1cm}|M{1cm}|}
			\hline
			Dataset
			& Algorithms & NOG (LR)
			&ACC (LR)  & NOG (BLR) &ACC (BLR)\\
			\hline
			\;  & {Sd-REG-LBFGS}  &\textbf{0.0288}   &\textbf{95.27\%} &\textbf{0.0294} &\textbf{95.47\%}\\
			\; & SdLBFGS       &0.0313   &95.17\% &0.3351 &91.67\%\\
			BNA & RSA           &0.0317   &95.11\% &0.3306 &90.36\%\\
			\;   & SAA           &1.7341   &95.11\% &1.6801 &90.35\%\\
			\;   & SGD           &0.0318   &95.10\% &0.3371 &90.35\%\\
			\;   & Adam          &0.2592   &92.23\% &0.1349 &94.52\%\\
			\hline
			
			\hline
			\;       & Sd-REG-LBFGS   &\textbf{0.010}   &97.11\% &\textbf{0.0073} &91.42\% \\
			\;     & SdLBFGS        &0.012   &97.31\% &0.0077 &91.43\%\\
			WINL & RSA       &0.0654   &95.89\% &0.023 &91.12\% \\
			\;   & SAA           &1.36    &95.89\% &0.595 &91.12\% \\
			\;   & SGD           &0.0653   &95.90\% &0.023 &91.12\%\\
			\;   & Adam          &0.54     &80.74\% &0.10 &86.46\%\\
			\hline

			\hline
			\;   & Sd-REG-LBFGS  &\textbf{0.013}   &\textbf{87.33\%} &\textbf{0.0188} &\textbf{88.21\%}\\
			\;   & SdLBFGS       &0.017   &86.98\% &0.0741 &86.82\%\\
			IONO & RSA           &0.087   &85.70\% &0.096 &84.59\%\\
			\;   & SAA           &0.795   &85.70\% &0.795 &84.87\%\\
			\;   & SGD           &0.087   &85.70\% &0.099 &84.87\%\\
			\;   & Adam          &0.19    &79.13\% &0.207 &80.80\%\\
			\hline
			
			\hline
			\;  & Sd-REG-LBFGS  &\textbf{0.017}   &87.55\% & \textbf{0.016}& 87.56\%\\
			\;& SdLBFGS       &0.020   &87.68\% &0.02 & 87.67\%\\
			ELEG  & RSA           &0.043   &86.72\% &0.024 &87.27\% \\
			\;   & SAA           &0.502   &86.72\% &0.502&87.27\%\\
			\;   & SGD           &0.042   &86.72\% &0.0241 &87.27\% \\
			\;   & Adam          &0.489   &52.69\% &0.493 &64.53\%\\
			\hline

		\end{tabular}
	\end{center}
\end{table}

\subsection{Numerical results using synthetic dataset}
In this subsection, we conduct the numerical experiments using synthetic dataset. For the binary classification schemes, we initialize the parameter to be optimized as  $\theta_0$, which is generated from a Gaussian distribution $\mathcal{N}(0,I)$. For SDVBI, the initial values of mean $\mu$ and the covariance matrix $S$ are set to  $\theta_0$ and identity matrix $S_0=I$,  respectively. We generate 5000 synthetic data points for 5-fold cross validation and 50 Monte Carlo runs in the following manner. Each of the sample $x_n$ is of dimension $d=50$ and is randomly drawn from a  uniform distribution  $[0,1]^{d}$. The desired parameter   $\bar{\theta}$ generated from the uniform distribution $[-1,1]^d$ is used to generate the true class labels $z_n=\mathbb{I}(\bar{\theta}^Tx_n>0)$ for each of the sample $x_n$. Using the synthetic dataset, we first consider the logistic regression problem and the objective function given in (\ref{logreg}). 

\subsubsection{Logistic Regression}
In Figs. 1(a) and 1(b), we illustrate the effect of batch size on the Sd-REG-LBFGS algorithm in terms of NOG and ACC, respectively. The regularized parameter $\gamma$ is set to $\gamma=10^{-4}$ and $\delta$ to $\delta=1.25\gamma+0.01$ correspondingly. Fig. 1(a) shows that the proposed approach consistently performs better than the SdLBFGS, SGD, RSA, SAA and Adam algorithms in NOG. Larger batch size generally leads to better performance for all algorithms. This is due to less variance of the stochastic gradient with larger batch size.

Figs. 1(a) and (b) show that the proposed approach and the SdLBFGS consistently performs better than SGD, RSA, SAA and Adam algorithms in terms of NOG and ACC, respectively.  Moreover, the proposed algorithm performs consistently well for different batch sizes, which suggests that the incorporation of regularization helps to reduce estimation variance and hence it is more robust to the variations of batch sizes.

In Figs. 2(a) and 2(b), we report the effect of various memory sizes on the performance of the proposed Sd-REG-LBFGS.  We set the step size constant to $r=7$ and the batch size $m=100$ for the proposed approach and SdLBFGS. The regularized parameters of the proposed approach are set to $\gamma=10^{-4}$ and $\delta$ to $\delta=1.25\gamma+0.01$, respectively,  which satisfies the condition $0.8\delta>\gamma$. Furthermore, the iteration interval length is set to $L=10$. From the figures, we can see that the proposed approach and the SdLBFGS give better NOG and ACC performance than the SGD, RSA, SAA and Adam. Moreover, a larger memory size generally lead to more accurate approximation of the Hessian matrix and hence a better performance.

In Figs. 3(a) and 3(b), we study the effect of the regularization parameter  on Sd-REG-LBFGS in terms of NOG and ACC, respectively.  The following values of $\gamma = 10^{-2},10^{-3},10^{-4}$ are employed.  We can see that the proposed approach performs better in terms of NOG and ACC. We notice the small amount of regularization imposed in the proposed Sd-REG-LBFGS method generally lead to better NOG than the SdLBFGS while its ACC is similar to SdLBFGS.

Overall, we find that the proposed approach performs better than other conventional algorithms in terms of NOG and ACC. This may be attributed by the small amount of regularization applied to the proposed approach, which improves the numerical stability and hence it converges closer to the stationary point (lower NOG). Meanwhile, we notice that the ACC of the proposed Sd-REG-LBFGS and the SdLBFGS algorithms are quite similar under this setting. We shall compare these algorithms more formally using a statistical test on their average classification accuracies at different settings in Section V-B. 

\captionsetup{font=footnotesize}

\begin{figure*} 
	\begin{minipage}[t]{0.325\linewidth} 
		\centering
		\begin{subfigure}{1\textwidth}
			\centering
			\includegraphics[width=1.0\linewidth]{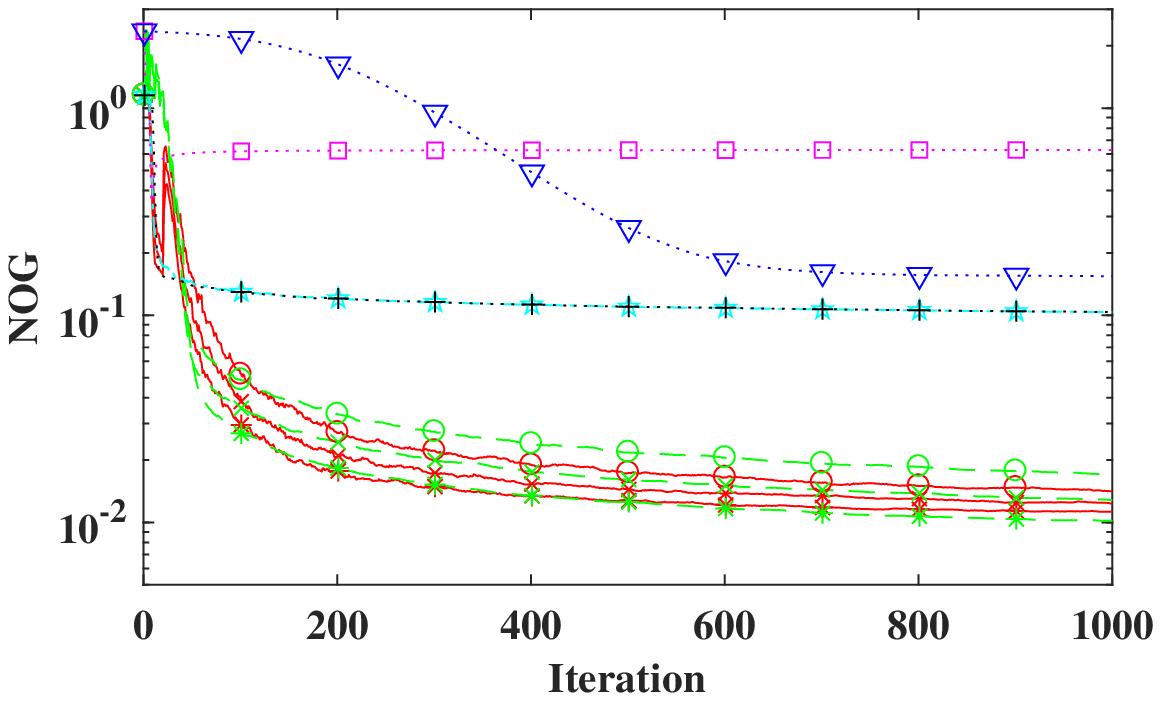}
			\label{fig:batch_logreg_syn_a}
		\end{subfigure} \\
		(a)
		\begin{subfigure}{1\textwidth}
			\centering
			\includegraphics[width=1.0\linewidth]{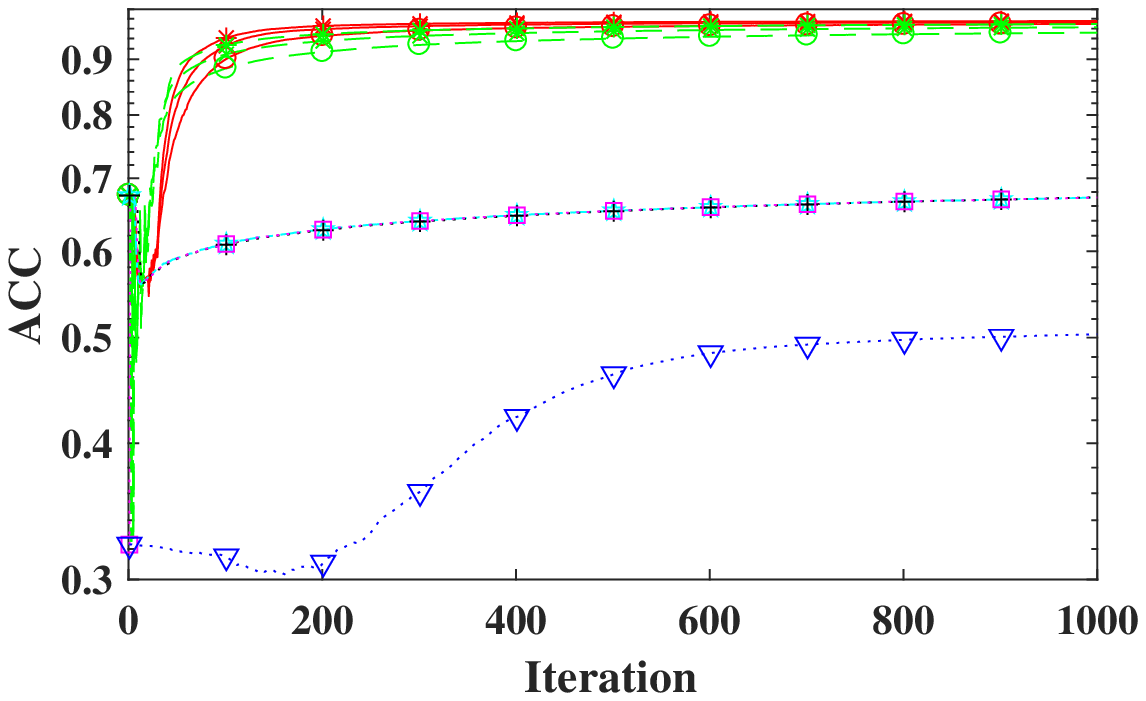}
			\label{fig:batch_logreg_syn_b}
		\end{subfigure} \\
		(b)
		\begin{subfigure}{1\textwidth}
			\centering
			\includegraphics[width=0.8\linewidth]{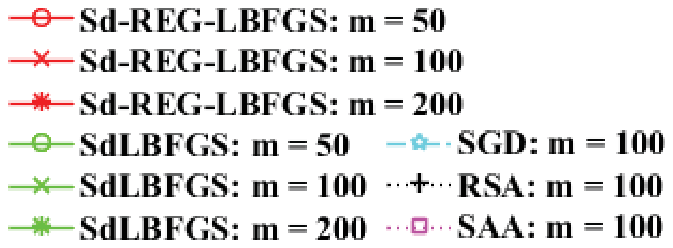}
			\label{}
		\end{subfigure}
		\caption{The (a) Norm of Gradient (NOG) and (b) Classification Accuracy (ACC) of logistic regression solved using various algorithms with different batch sizes averaged over 50 Monte Carlo simulations. The synthetic dataset is used.}
		\label{fig:batch_logreg_syn}
	\end{minipage} \hfill
	\begin{minipage}[t]{0.325\linewidth} 
		\centering
		\begin{subfigure}{1\textwidth}
			\centering
			\includegraphics[width=1.0\linewidth]{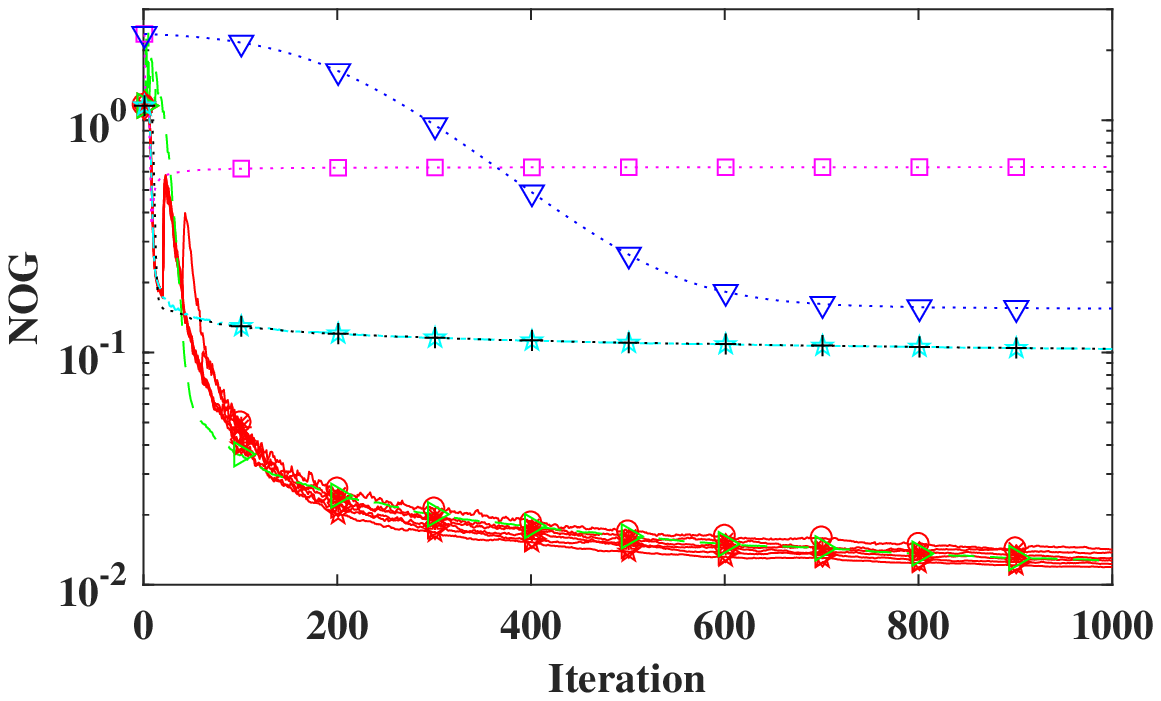}
			\label{fig:mem_logreg_syn_a}
		\end{subfigure}\\
		(a)
		\begin{subfigure}{1\textwidth}
			\centering
			\includegraphics[width=1.0\linewidth]{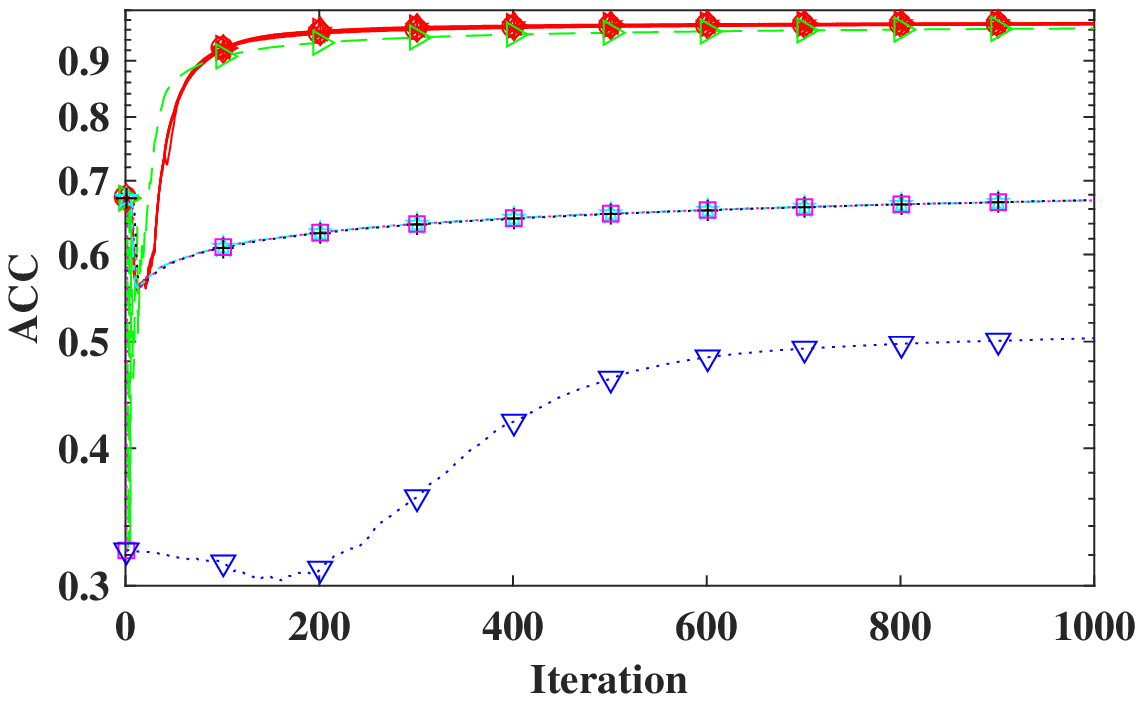}
			\label{fig:mem_logreg_syn_b}
		\end{subfigure}\\
		(b)
		\begin{subfigure}{1\textwidth}
			\centering
			\includegraphics[width=0.9\linewidth]{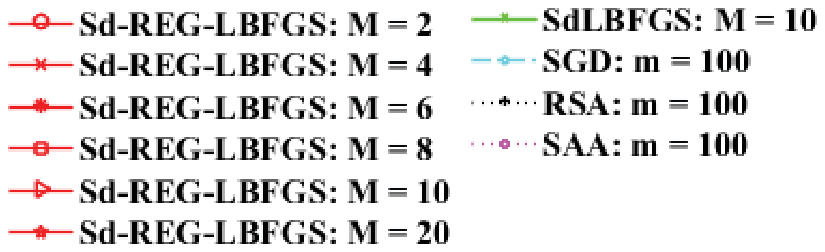}
			\label{}
		\end{subfigure}
		\caption{The (a) NOG and (b) ACC of logistic regression solved using various algorithms with different memory sizes averaged over 50 Monte Carlo simulations.  For comparison, SdLBFGS ,  SGD, RSA, SAA and Adam are implemented. The synthetic dataset is used.}
		\label{fig:mem_logreg_syn}
	\end{minipage} \hfill 
	\begin{minipage}[t]{0.325\linewidth} 
		\centering
		\begin{subfigure}{1\textwidth}
			\centering
			\includegraphics[width=1.0\linewidth]{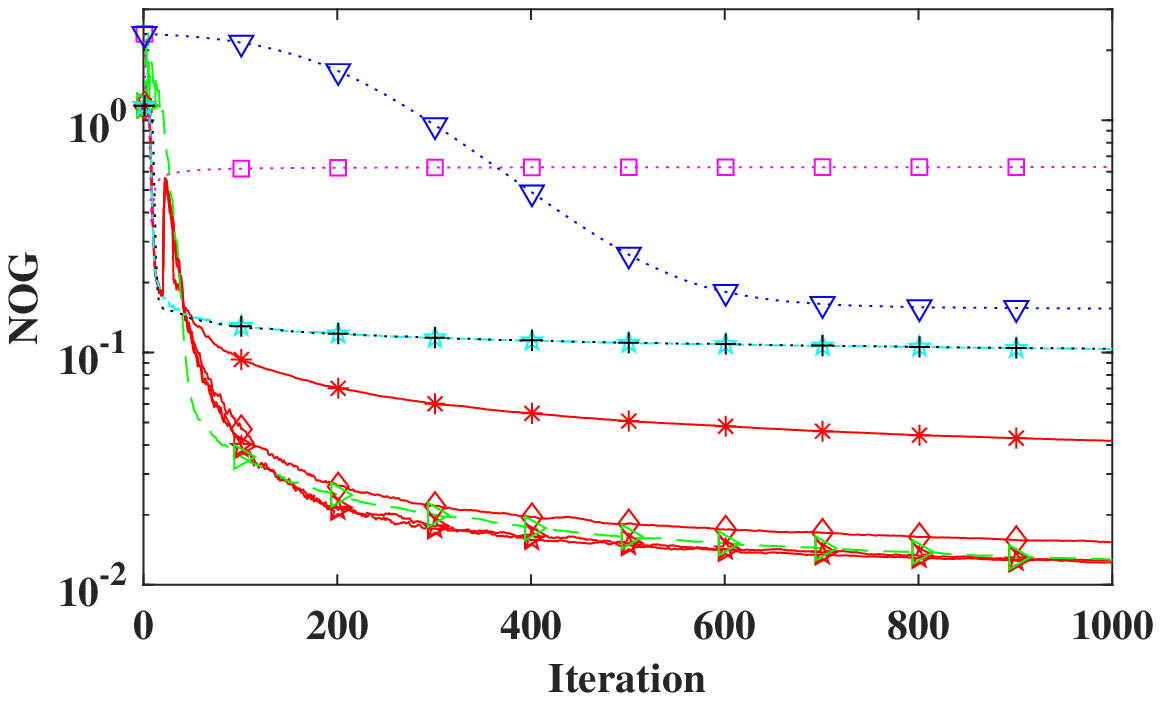}
			\label{fig:gamma_logreg_syn_a}
		\end{subfigure}\\
		(a)
		\begin{subfigure}{1\textwidth}
			\centering
			\includegraphics[width=1.0\linewidth]{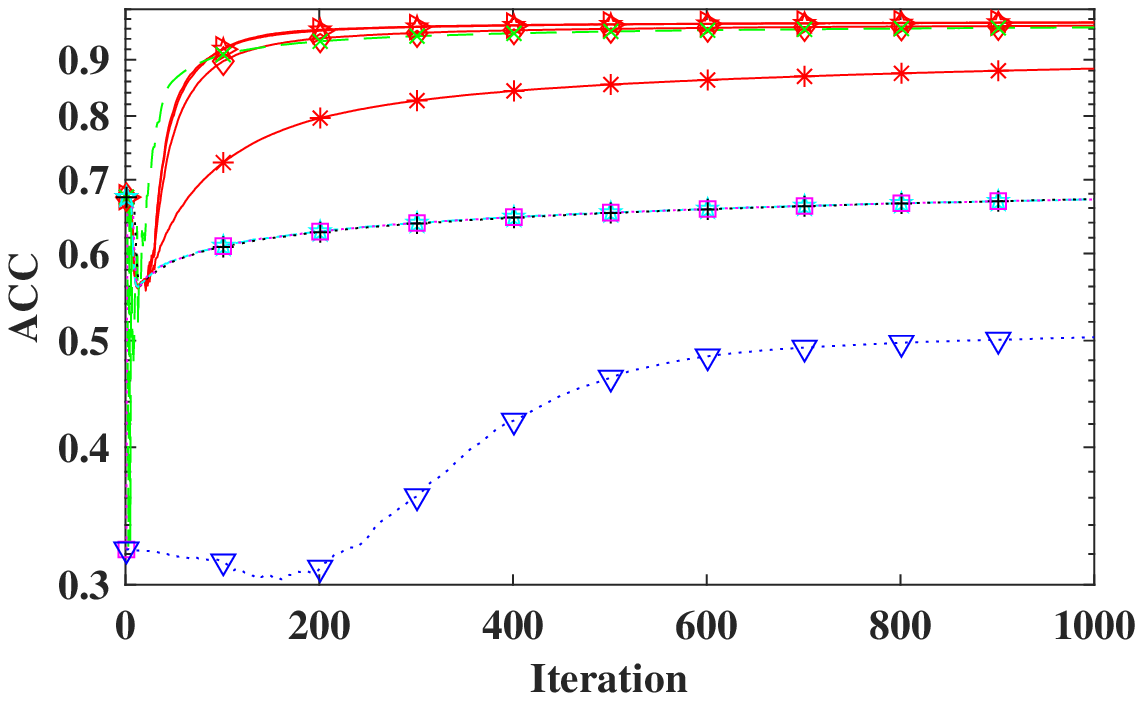}
			\label{fig:gamma_logreg_syn_b}
		\end{subfigure}\\
		(b)
		\begin{subfigure}{1\textwidth}
			\centering
			\includegraphics[width=1\linewidth]{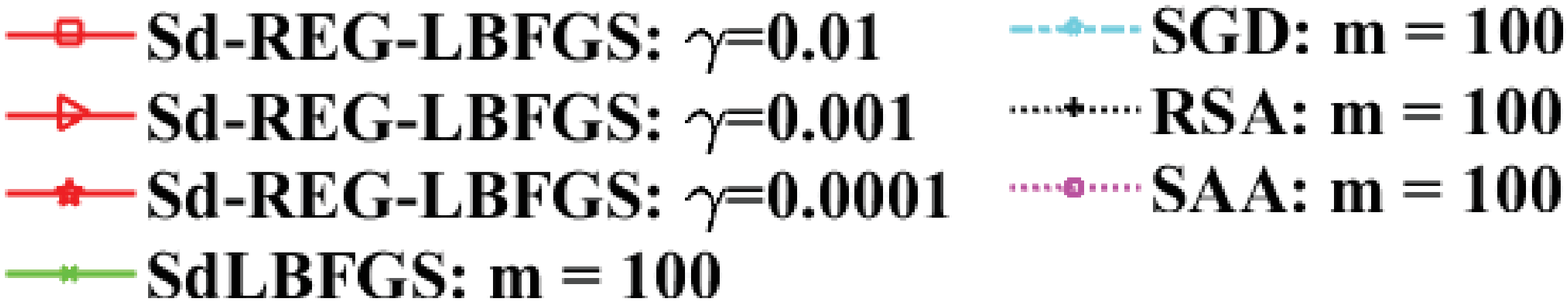}
			\label{}
		\end{subfigure}
		\label{fig:gamma_logreg_syn}
		\caption{The effect of regularized parameter $\gamma$ on the (a) NOG and (b) ACC of logistic regression solved using the proposed Sd-REG-LBFGS. For comparison, SdLBFGS ,  SGD, RSA, SAA and Adam are implemented. The synthetic dataset is used.}
	\end{minipage}%
	
\end{figure*}

\subsubsection{SDVBI for Bayesian Logistic Regression}

In this subsection, numerical experiments are performed on SDVBI for Bayesian logistic regression using synthetic dataset. We will study various values of the batch size $m$, the memory size $M$ and the regularization parameter $\gamma$, under which Sd-REG-LBFGS is performed for the optimization. 

Figs 4(a) and 4(b) show the NOG and ACC of SDVBI, respectively, solved using different algorithms with various batch sizes. In this experiment, we fix the regularization parameters to  $\gamma=10^{-4}$ and   $\delta=1.25\gamma+0.01$ respectively. Moreover, we set the memory size to  $M=10$ and the interval length to  $L=10$. 

From the figures, it can be seen that the proposed method generally outperforms the Sd-LBFGS, SGD, RSA, SAA and Adam algorithms with all batch sizes studied. Moreover, the proposed algorithm performs consistently well for different batch size, which suggests the incorporation of regularization helps to reduces estimation variance and hence it is more robust to the variations of batch sizes. This enables us to choose a smaller batch size so that it could reduce computational cost without sacrificing much classification performance of the SDVBI in Bayesian logistic regression.

 Figs. 5(a) and 5(b) show the effect of memory size on the NOG and classification performances of SDVBI in Bayesian logistic regression using the proposed Sd-REG-LBFGS. The Sd-LBFGS, SGD, RSA, SAA and Adam are also included as benchmarks. Similar to previous sub-sections, we fix the regularized parameter to  $\gamma=10^{-4}$ and   $\delta=1.25\gamma+0.01$ respectively. Moreover, the batch size is set to $m=100$  and the interval length for Sd-REG-LBFGS is set to  $L=10$. 

From the figures, we find that the NOG and ACC performance of the proposed approach is generally better than other approaches.  Thus, we can choose a relatively small memory size to reduce the computational cost without sacrificing performance.

In Figs 6(a) and 6(b), we report the effect of regularization parameter $\gamma$  on Sd-REG-LBFGS for SDVBI. In general, smaller $\gamma$  value yields better performance in terms of NOG. For the Sd-REG-LBFGS with  $\gamma=10^{-2},\:10^{-3}\;10^{-4}$, the improvement in classification performance decreases when   decreases. We notice the small amount of regularization imposed in the proposed Sd-REG-LBFGS method generally lead to better NOG than the SdLBFGS while its ACC is similar to SdLBFGS. 

Overall, the proposed Sd-REG-LBFGS performs better than the SdLBFGS, SGD, RSA, SAA and Adam algorithms in terms of NOG and classification accuracy. Moreover, the classification performance of the proposed Sd-REG-LBFGS algorithm is less vulnerable to insufficient samples caused by small batch size as regularization is imposed to avoid ill-conditioning of the Hessian update. On the other hand, the proposed Sd-REG-LBFGS and SdLBFGS gives similar performance when the number of samples is large.

Regarding to the choice of the algorithmic parameters including the batch size  $m$, memory size $M$  and the regularization parameter  $\gamma$, we observe a choice of $m=100$  yields the best performance for most algorithms under the datasets we have considered. For the proposed Sd-REG-LBFGS method and the SdLBFGS method, a memory size of $M=8$ will suffice. Beyond these values, the performance improvement is not so significant. Moreover, the complexity and computational time increases with the two parameters and hence it is desirable to keep them as small as possible. A small amount of regularization, such as  $\gamma=10^{-4}$, is adequate to reduce the fluctuations under sufficient samples.

We notice that the proposed Sd-REG-LBFGS and SdLBFGS algorithms gives similar performance with sufficient large number of samples. We shall further compare the two algorithms more rigorously using a statistical test under different settings in Section V-C.

\begin{figure*} 
	\begin{minipage}[t]{0.325\linewidth} 
		\centering
		\begin{subfigure}{1\textwidth}
			\centering
			\includegraphics[width=1.0\linewidth]{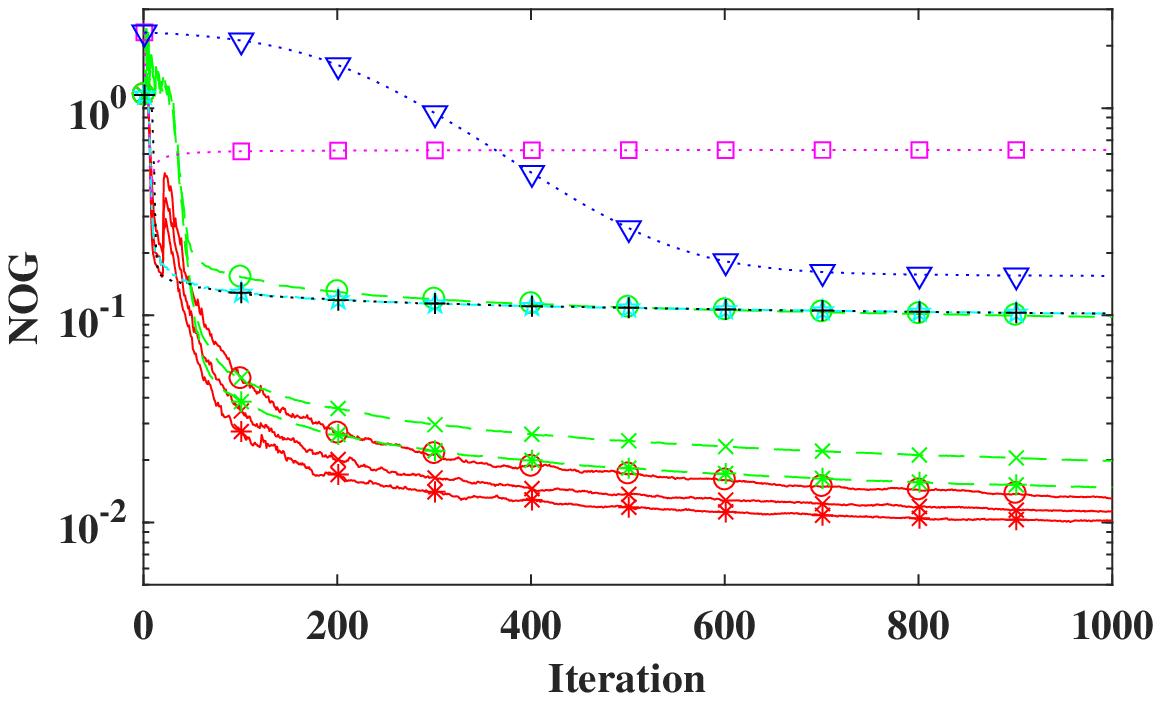}
			\label{fig:batch_bayeslogreg_syn_a}
		\end{subfigure}\\
		(a)
		\begin{subfigure}{1\textwidth}
			\centering
			\includegraphics[width=1.0\linewidth]{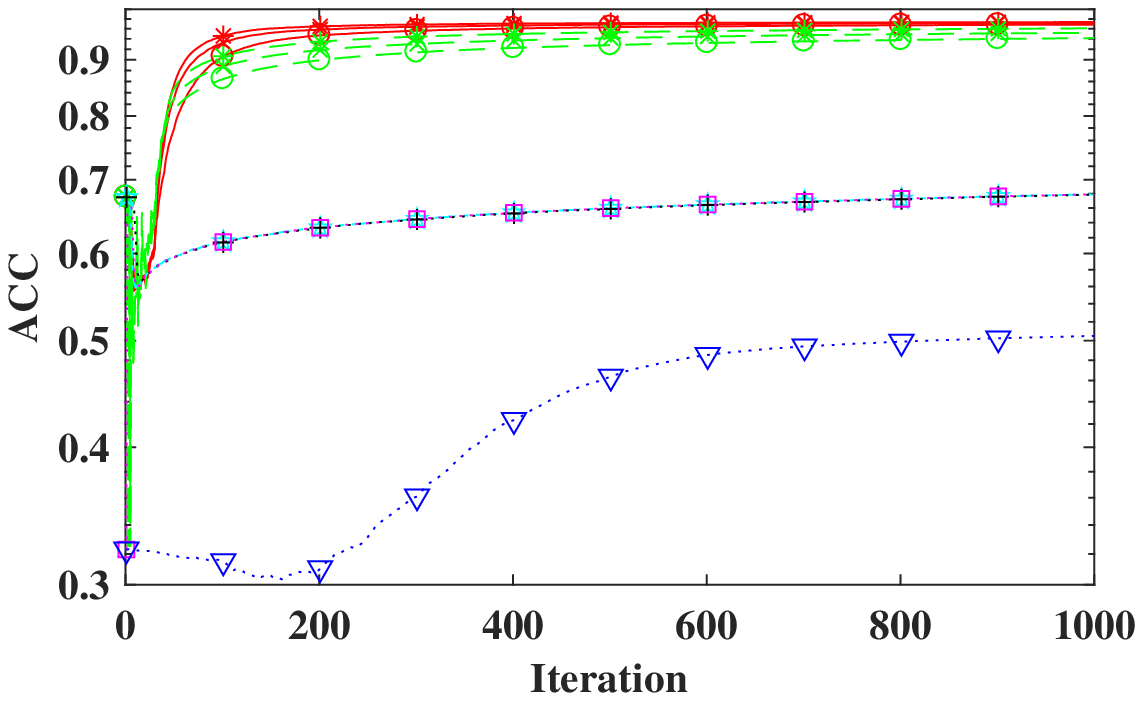}
			\label{fig:batch_bayeslogreg_syn_b}
		\end{subfigure}\\
		(b)
		\begin{subfigure}{1\textwidth}
			\centering
			\includegraphics[width=0.8\linewidth]{Fig1_legend.eps}
			\label{fig:batch_bayeslogreg_syn_c}
		\end{subfigure}
		\caption{The (a) NOG and (b) ACC of different algorithms with various batch sizes in solving SDVBI in Bayesian logistic regression. For comparison, SdLBFGS, SGD, RSA, SAA and Adam are included. The synthetic dataset is used.}
		\label{fig:batch_bayeslogreg_syn}
	\end{minipage} \hfill 
	\begin{minipage}[t]{0.325\linewidth} 
		\centering
		\begin{subfigure}{1\textwidth}
			\centering
			\includegraphics[width=1.0\linewidth]{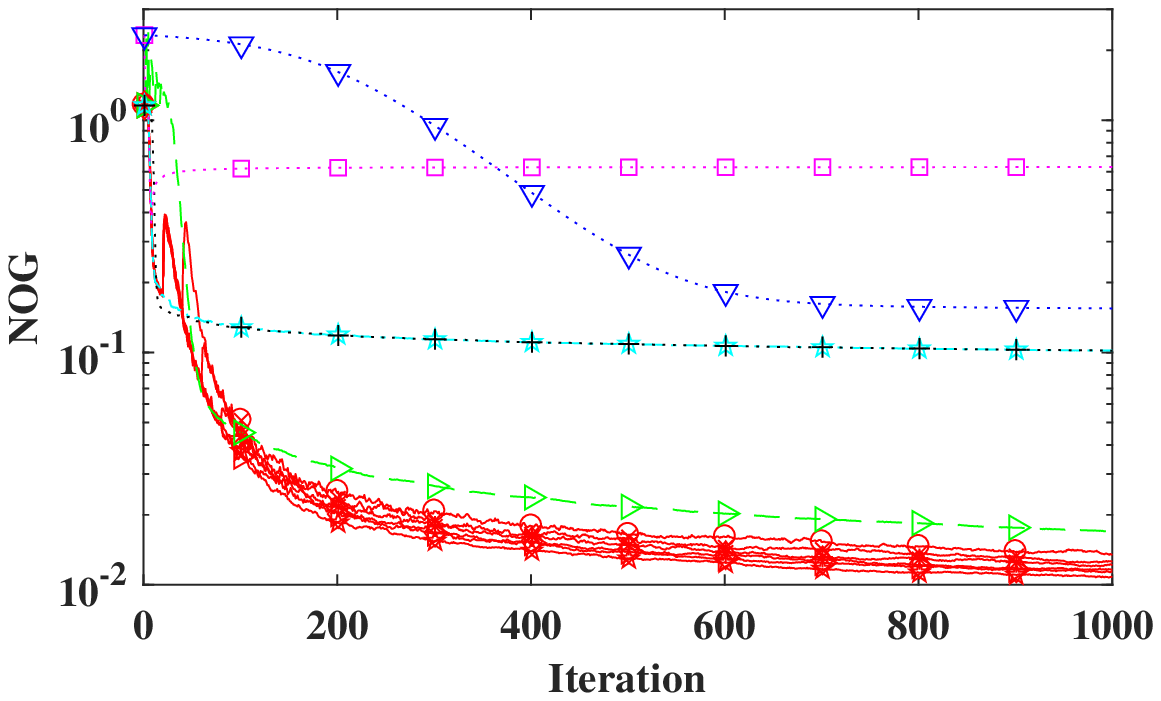}
			\label{fig:mem_bayeslogreg_syn_a}
		\end{subfigure}\\
		(a)
		\begin{subfigure}{1\textwidth}
			\centering
			\includegraphics[width=1.0\linewidth]{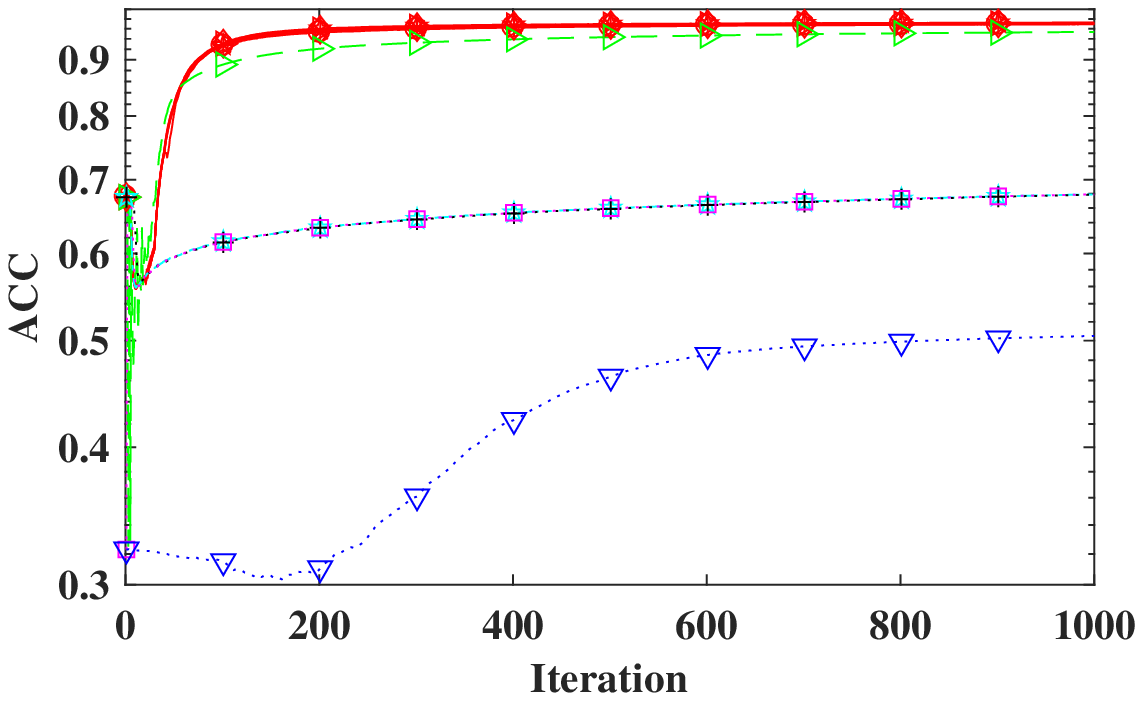}
			\label{fig:mem_bayeslogreg_syn_b}
		\end{subfigure}\\
		(b)
		\begin{subfigure}{1\textwidth}
			\centering
			\includegraphics[width=0.9\linewidth]{Fig2_legend.eps}
			\label{fig:mem_bayeslogreg_syn_c}
		\end{subfigure}
		\caption{The (a) NOG and (b) ACC of SDVBI in Bayesian logistic regression solved using various algorithms with different memory sizes averaged over 50 Monte Carlo simulations.  The synthetic dataset is used.}
		\label{fig:mem_bayeslogreg_syn}
	\end{minipage} \hfill
	\begin{minipage}[t]{0.325\linewidth} 
		\centering
		\begin{subfigure}{1\textwidth}
			\centering
			\includegraphics[width=1.0\linewidth]{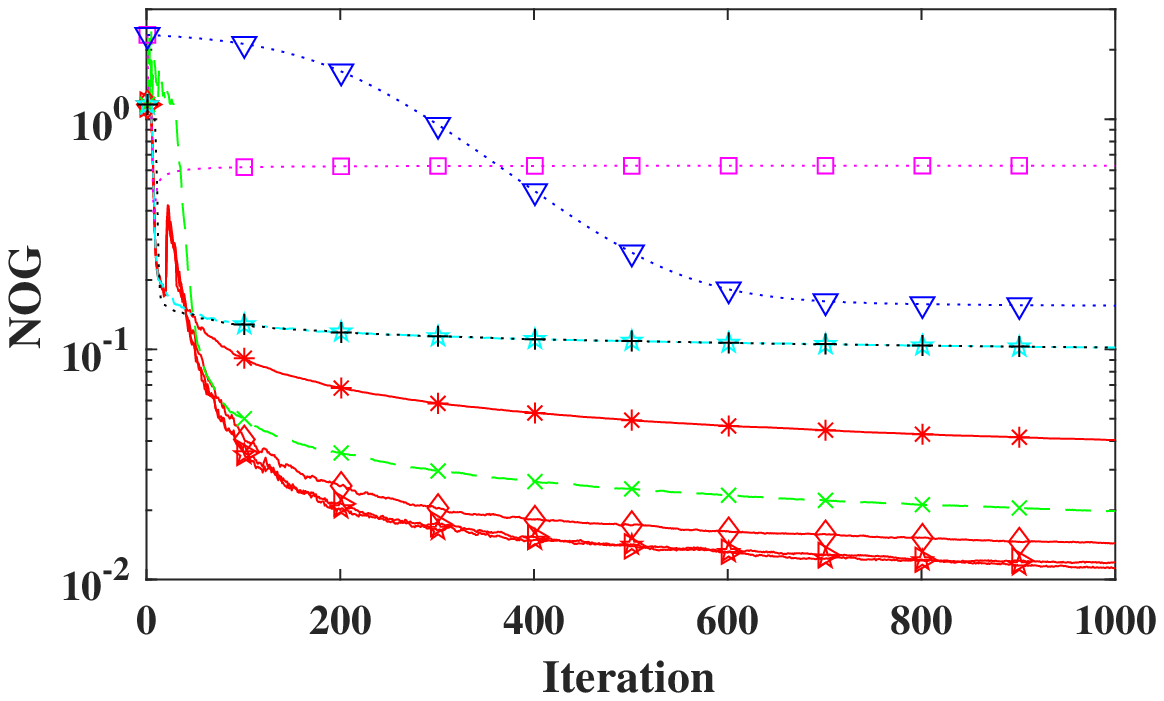}
			\label{fig:gamma_bayeslogreg_syn_a}
		\end{subfigure}\\
		(a)
		\begin{subfigure}{1\textwidth}
			\centering
			\includegraphics[width=1.0\linewidth]{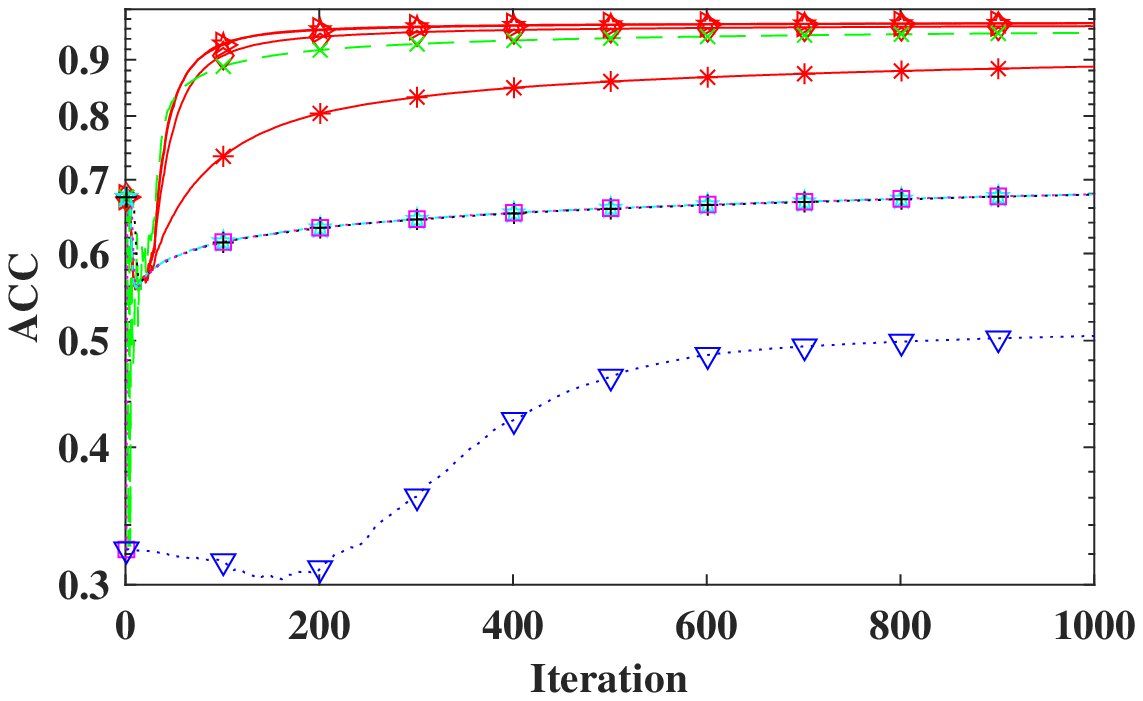}
			\label{fig:gamma_bayeslogreg_syn_b}
		\end{subfigure}\\
		(b)
		\begin{subfigure}{1\textwidth}
			\centering
			\includegraphics[width=1\linewidth]{Fig3_legend.eps}
			\label{fig:gamma_bayeslogreg_syn_c}
		\end{subfigure}
		\caption{The effect of regularized parameter $\gamma$ on the (a) NOG and (b) ACC of SDVBI in Bayesian logistic regression solved using the proposed Sd-REG-LBFGS. For comparison, SdLBFGS, SGD, RSA, SAA and Adam are implemented. The synthetic dataset is used.}
		\label{fig:gamma_bayeslogreg_syn}
	\end{minipage}%
\end{figure*}

\subsection{Comparison of classification performance of various algorithms using Statistical Significance Testing}

In this section, we employ nonparametric statistical tests - \textit{sign test} and \textit{Wilcoxon  paired-difference test} for the evaluation of the statistical significance of whether the proposed algorithm performs significantly better than the SGD, RSA, SAA and Adam algorithms on average, or vice versa, in terms of classification accuracy. It should be noted that the $t$-test may not be a proper choice as classification accuracies (ACC) are bounded and hence they are not normally distributed \cite{R48,R49}. First, Table \ref{avg_acc} shows the average classification accuracy for each algorithm over all batch sizes, and the parameters are set for each algorithm as follows: 

1. Batch size:  $m=5,\;10,\;30,\;50,\;100,\;200,\;$ for synthetic dataset scenario, and  $m=5,\;10,\;20,\;30,\;50,\;100,\;$ for \textit{scene} dataset scenario; 

2. Memory size: $M=10$ for Sd-REG-LBFGS and SdLBFGS;

3. Step size constant $r$: $r=7$ for all algorithms;

4. Regularization parameters:  $\gamma=10^{-4}$ and $\delta=1.25\gamma+0.01$.

Here, we abbreviate Sd-REG-LBFGS and SdLBFGS as SRL and SDL for convenience respectively. More precisely, the \textit{sign} test tests the following hypotheses: 

\begin{equation}
	H_0:\;\mu_X-\mu_Y=0\;\;\text{vs}\;\;H_1:\;\mu_X-\mu_Y>0,
\end{equation}
where $\mu_X$ and $\mu_Y$ are the median classification accuracies of algorithms A and B, respectively. The test statistic of the sign test is given as 
\begin{equation}
	T_S:\;\text{number of times that }x_i-y_i>0,
\end{equation}
where $x_i$ and $y_i$ are the classification accuracies of algorithms A and B for the \textit{i}-th experiment, respectively. The one-sided \textit{p}-value can be obtained by a binomial test as $P=\text{Pr}(T_S\geq t_S|H_0)=\sum_{i=t_S}^{n} \binom{n}{i}0.5^n$, where $t_S$ is the observed number of times that $x_i-y_i>0$. $n$ is the total number of experiments performed.

For \textit{Wilcoxon paired-difference} test, the following hypotheses are considered
\begin{equation}
	H_0:\;|x_i-y_i|\text{ follows a symmetric distribution around zero,}
\end{equation}
\begin{equation}
\begin{aligned}
H_1:&\;|x_i-y_i|\text{ does not follows a symmetric distribution aro-}\\
&\text{ und zero.}
\end{aligned}
\end{equation}
The test statistic is given as $T_W=\sum_{i=1}^{n_R}\text{sign}(x_i> y_i) R_i$, where $\text{sign}(x> y)$ is defined to be
\begin{equation}\label{signfunc}
\text{sign}(x> y)=\left\{
\begin{aligned}
&+1,\text{ if} (x> y)\\
&-1, \;otherwise,
\end{aligned}
\right.
\end{equation}
and $R_i$ is the rank order of $|x_i-y_i |$.  $n_R$ is the number of experiments after excluding those with $|x_i-y_i |=0$. For $n_R<20$, the exact distribution is used. For $n_R\geq20$, a $z$-score can be calculated as $z=T_w/\sigma_W$, where $\sigma_W=\sqrt{n_R (n_R+1)(2n_R+1)/6}$. The right-sided $p$-value for $x_i>y_i$ is $P=\text{Pr}(T_W\geq t_w|H_0)$, where $t_w$ is the observed sum of rank. The $p$-values are obtained using MATLAB function \textit{signrank}. The batch size is set to $m=5$ as our proposed method is robust and efficient in particular for small batch sizes. 
 
The results of \textit{sign test} are shown in Table \ref{sign_p}. The $\text{log}\;p$-values for \textit{Wilcoxon paired-difference test} are shown in Table \ref{wilcoxon}.  The batch size is set to $m=5$ as our proposed method is robust and efficient in particular for small batch sizes. From the tables, we can see that the proposed approach obtains the highest ACC with statistical significance and the mean difference in ACC between the proposed approach and other algorithms is statistically significant for $\text{log }p<-1.3$, (a.k.a.  $p<0.05$). A key observation is that we find that the proposed approach performs much better than the SdLBFGS under small batch size. This is possibly attributed to the incorporation of the proposed regularization scheme, which is useful to improve numerical stability under small sample size. For sufficient samples, the performance of our algorithm is similar to SdLBFGS. Such observations can be found in the sensitivity study of the different parameters, which is omitted here due to page limitation. Interested readers are referred to Section III of the supplementary material for details.

\begin{table}[!t]
\renewcommand{\arraystretch}{1.3}
\caption{The average classification accuracy of each algorithm over all batch sizes.}
\label{avg_acc}
\centering

\begin{tabular}{|M{0.6cm}|M{0.83cm}|M{0.83cm}|M{0.83cm}|M{0.83cm}|M{0.83cm}|M{0.83cm}|}
\hline
ACC & SRL &SDL &SGD &RSA &SAA &Adam \\
\hline
S1 & 95.14\% &89.67\% &66.69\% &66.58\% & 66.70\%&50.25\%\\
\hline
S2 & 95.25\% & 89.78\% & 67.26\% & 67.23\% &67.26\% &50.35\%\\
\hline
S3 & 80.80\% & 80.48\%&77.32\% &77.36\% &77.32\% &65.72\%\\
\hline
S4 & 80.90\% &79.12\% &76.38\%& 76.50\%&76.38\% &66.27\%\\
\hline
\end{tabular}
\end{table}

%

\begin{table}[!t]
	\renewcommand{\arraystretch}{1.3}
	\caption{ Right-sided $\text{log }p$ values obtained from \textit{sign test} on mean classification accuracy of various algorithms averaged over 50 Monte Carlo simulations}
	\label{sign_p}
	\centering
	
	\begin{tabular}{|M{0.9cm}|M{0.9cm}|M{0.9cm}|M{0.9cm}|M{0.9cm}|M{0.9cm}|}
		\hline
		$\text{log }p$
		& SRL vs
		SDL
		&SRL vs SGD &SRL vs RSA &SRL vs SAA &SRL vs Adam\\
		\hline
		S1 & -15.05 &-15.05 &-15.05 & -15.05& -15.05\\
		\hline
		S2 &  -15.05 & -15.05 & -15.05 &-15.05 &-15.05\\
		\hline
		S3 &  -15.05 &-15.05 &-15.05 &-15.05 &-15.05 \\
		\hline
		S4 & -10.732 &-15.05&-15.05&-15.05&-15.05 \\
		\hline
	\end{tabular}
\end{table}

\begin{table}[!t]
	\renewcommand{\arraystretch}{1.3}
	\caption{\textit{Wilcoxon paired-difference test} on mean classification accuracy of various algorithms averaged over 50 Monte Carlo simulations}
	\label{wilcoxon}
	\centering
	
	\begin{tabular}{|M{0.6cm}|M{1.1cm}|M{0.9cm}|M{0.9cm}|M{0.9cm}|M{0.9cm}|}
		\hline
		$\text{log}\;p$
		& SRL vs
		SDL
		&SRL vs SGD &SRL vs RSA &SRL vs SAA &SRL vs Adam\\
		\hline
		S1 & -9.40 &-9.40 &-9.40 & -9.40& -9.40\\
		\hline
		S2 &  -9.40 & -9.41 & -9.41 &-9.40 &-9.40\\
		\hline
		S3 &  -9.41 &-9.41 &-9.41 &-9.41 &-9.41 \\
		\hline
		S4 & -9.02 &-9.41 &-9.42 &-9.41&-9.42 \\
		\hline
	\end{tabular}
\end{table}

\section{Conclusion}

A novel Sd-REG-LBGS method for solving nonconvex and ill-conditioned stochastic optimization problems has been presented.  The convergence of the proposed method is established under reasonable assumptions. The effectiveness of the proposed method is studied via the logistic regression and Bayesian logistic regression problems in machine learning for both synthetic and real datasets. The effect of using different algorithmic parameters is also studied.  Experimental results show that the proposed Sd-REG-LBFGS method generally outperforms SdLBFGS and exhibits superior performance for problems with small sample sizes. Moreover, the proposed method is less sensitive to the variations of the batch size and memory size than the SdLBFGS method. For future work, we shall consider the extension of our method to distributed optimization \cite{MokhtariAryan2014RRSB2,scchan1,scchan2,scchan3} and asynchronous distributed optimization \cite{AsynADMM,AsynBlockADMM}.


%

%
%
%
%
%

\ifCLASSOPTIONcaptionsoff
  \newpage
\fi


\begin{thebibliography}{1}
	%
	
	
	
	\bibitem{SBBFGS}
	R. M. Gower, D. Goldfarb, and P. Richtarik, \textquotedblleft{}Stochastic block BFGS: Squeezing more curvature out of data,\textquotedblright in \textit{33rd Proc. Int. Conf. Mach. Learn}, 1869-1878, June 19-24, 2016.

	\bibitem{SGD-QN}
	A. Bordes, L. Bottou, and P. Gallinari, \textquotedblleft{}SGD-QN: Careful Quasi-Newton Stochastic Gradient Descent,\textquotedblright \textit{J. Mach. Learn. Res.}, 10,  pp. 1737 – 1754, Jul. 2009.
	
	\bibitem{MokhtariAryan2014RRSB}
	A. Mokhtari and A. Ribeiro, \textquotedblleft{}RES: Regularized Stochastic BFGS Algorithm,\textquotedblright \textit{IEEE Trans. Signal Process.}, vol. 62, no. 23,  pp. 6089 - 6104, Dec.1, 2014.
	
	\bibitem{MokhtariAryan2014RRSB2}
	M. Eisen, A. Mokhtari and A. Ribeiro, \textquotedblleft{}Decentralized Quasi-Newton Methods,\textquotedblright \textit{IEEE Trans. Signal Process.}, vol. 65, no. 10,  pp. 2613 - 2628, May, 2017.
	
	\bibitem{DSOCS}
    M. Neely, \textquotedblleft{}Distributed Stochastic Optimization
	via Correlated Scheduling,\textquotedblright \textit{IEEE/ACM Trans. Netw.}, vol. 24, no. 2,  pp. 759 - 772, April, 2016.
	
	
	\bibitem{PSI}
	P. Si, J. Yang, S. Chen, and H. Xi, \textquotedblleft{}Smoothness Constraint Based Stochastic Optimization
	for Wireless Scalable Video Streaming,\textquotedblright \textit{IEEE Commun.
	Lett.}, vol. 19, no. 5,  pp. 759 - 762, May, 2015.
	
    \bibitem{ARibeiro}
	A. Ribeiro, \textquotedblleft{}Ergodic Stochastic Optimization Algorithms for Wireless Communication and Networking,\textquotedblright \textit{IEEE Trans. Signal Process.}, vol. 58, no. 12,  pp. 6369 - 6386, Dec, 2010.
	
	
	\bibitem{SQN}
	X. Wang, S. Ma, D. Goldfarb, and W. Liu, \textquotedblleft{}Stochastic Quasi-Newton Methods for Nonconvex Stochastic Optimization,\textquotedblright \textit{SIAM J. Optim.}, vol. 27, no. 2,  pp. 927 – 956, 2017.
	
	\bibitem{IQN}
	A. Mokhtari, M. Eisen and A. Ribeiro, \textquotedblleft{}IQN: An Incremental Quasi-Newton Method with Local Superlinear Convergence Rate,\textquotedblright \textit{SIAM J. Optim.}, vol. 28, no. 2,  pp. 1670 – 1698, 2018.
	
	\bibitem{RHBYRD}
	R. H. Byrd, S. L. Hansen, J. Nocedal, and Y. Singer, \textquotedblleft{}A Stochastic Quasi-Newton Method for Large-Scale Optimization,\textquotedblright \textit{SIAM J. Optim.}, vol. 26, no. 2,  pp. 1008 - 1031, 2016.
	
	\bibitem{RHBYRD2}
	R. H. Byrd, G. M. Chin, W. Neveitt, and J. Nocedal, \textquotedblleft{}On the Use of Stochastic Hessian Information in Optimization Methods for Machine Learning
	,\textquotedblright \textit{SIAM J. Optim.}, vol. 21, no. 3,  pp. 977 – 995, Jan. 2011.
	
	\bibitem{OPREVIEW}
	L. Bottou, F. E. Curtis and J. Nocedal, \textquotedblleft{}Optimization Methods for Large-Scale Machine Learning,\textquotedblright \textit{SIAM Rev.}, vol. 60, no. 2,  pp. 223 - 311, 2018.
	
	
	
	\bibitem{Powell1978}
	M. J. D. Powell, \textquotedblleft{}Algorithms for nonlinear constraints that use lagrangian functions,\textquotedblright \textit{Math. Programming}, vol. 14, no. 1,  pp. 224 – 248, Dec., 1978.
	
	\bibitem{nocedal2006numerical}
	J. Nocedal, S. J. Wright,  \textit{Numerical Optimization}, New York:Springer-Verlag,  1999.
	
	\bibitem{oLBFGS}
	A. Mokhtari, and A. Ribeiro, \textquotedblleft{}Global Convergence of Online Limited Memory BFGS,\textquotedblright \textit{J. Mach. Learn. Res.},   vol. 16, no. 1,  pp. 3151 - 3181, Jan., 2015.
	\bibitem{oLBFGS2}
	N. Schraudolph, J. Yu, and S. Gunter, \textquotedblleft{}A stochastic quasi-Newton method for online convex optimization,\textquotedblright in \textit{Proc. 11th Int. Conf. Artif. Intell. Statist.},  pp. 433 –	440, 2007.
	
	\bibitem{Bishop}
	C. Bishop, \textit{Pattern Recognition and Machine Learning},   Springer New York., 2006.
	
	\bibitem{NONCON}
	C. Wang, and D. M. Blei , \textquotedblleft{}Variational Inference in Nonconjugate Models,\textquotedblright \textit{J. Mach. Learn. Res.},   vol. 14, no. 1,  pp. 1005 - 1031, Jan., 2013.
	
	\bibitem{StoFirZeroSP}
	S. Ghadimi and G. Lan, \textquotedblleft{}Stochastic First- and Zeroth-Order Methods for Nonconvex Stochastic Programming,\textquotedblright \textit{SIAM J. Optim.}, vol. 23, no. 4,  pp. 2341 - 2368, 2013.
	
	\bibitem{RobustSA}
	A. Nemirovski and A. Juditsky and G. Lan and A. Shapiro, \textquotedblleft{}Robust Stochastic Approximation Approach to Stochastic Programming,\textquotedblright \textit{SIAM J. Optim.}, vol. 19, no. 4,  pp. 1574 - 1609, 2009.
	
	\bibitem{SBMD}
	C. Dang and G. Lan, \textquotedblleft{}Stochastic Block Mirror Descent Methods for Nonsmooth and Stochastic Optimization,\textquotedblright \textit{SIAM J. Optim.}, vol. 25, no. 2,  pp. 856 - 881, 2015.
	
	\bibitem{RobbinsHerbert1951ASAM}
	H. Robbins and S. Monro, \textquotedblleft{}A Stochastic Approximation Method,\textquotedblright \textit{Ann. Math. Statist.}, vol. 22, no. 3,  pp. 400 - 407, 1951.
	

	\bibitem{HoffmanSVI}
	M. D. Hoffman, D. M. Blei, C. Wang and J. Paisley, \textquotedblleft{}Stochastic Variational Inference,\textquotedblright \textit{J. Mach. Learn. Res.},   vol. 14, no. 1,  pp. 1303 - 1347, Jan., 2013.
	
	\bibitem{ZoubinGhahramani2015Pmla}
	Z. Ghahramani, \textquotedblleft{}Probabilistic machine learning and artificial intelligence,\textquotedblright \textit{Nature},   vol. 521, no. 7553,  pp. 452 – 459, 2013.
	
	\bibitem{BleiDavidM.2017VIAR}
	D. M. Blei, A. Kucukelbir and J. D. McAuliffe, \textquotedblleft{}Variational Inference: A Review for Statisticians,\textquotedblright \textit{J. Am. Statist. Assoc.},   vol. 112, no. 518,  pp. 859-877, 2017.
	
	\bibitem{scchan1}
	L. Zhang, H. C. Wu, C. H. Ho, S. C. Chan, \textquotedblleft{}A Multi-Laplacian Prior and Augmented Lagrangian Approach to the Exploratory Analysis of Time-Varying Gene and Transcriptional Regulatory Networks for Gene Microarray Data
	, to appear in \textquotedblright \textit{IEEE/ACM Trans. Comput. Biol. Bioinf.}. 
	
	\bibitem{scchan2}
	S. C. Chan, L. Zhang, H. C. Wu, and K. M. Tsui, \textquotedblleft{}A maximum
	a posteriori probability and time-varying approach for inferring
	gene regulatory networks from time course gene microarray
	data,\textquotedblright \textit{IEEE/ACM Trans. Comput. Biol. Bioinf.}, vol. 12, no. 1, pp. 123–135, 2015.
	
	\bibitem{scchan3}
	S. C. Chan, H. C. Wu, C. H. Ho and L. Zhang, \textquotedblleft{}An Augmented Lagrangian Approach for Distributed Robust Estimation in Large-Scale Systems, to appear in \textquotedblright \textit{IEEE Systems Journal}.
	
	\bibitem{VBISS}
	J. Paisley, D. M. Blei and M. I. Jordan, \textquotedblleft{}Variational Bayesian Inference with Stochastic Search,\textquotedblright in \textit{29th Proc. Int. Conf. Mach. Learn}, vol. 14,  pp. 1363 - 1370, 2012.
	
	\bibitem{YiSun2009Ssut}
	Sun Yi,	D. Wierstra, T. Schaul and J. Schmidhuber, \textquotedblleft{}Stochastic search using the natural gradient,\textquotedblright in \textit{26th Proc. Int. Conf. Mach. Learn},  vol. 382, pp. 1161 - 1168, 2009.
	
	\bibitem{NIPS2013_4937}
	R. Johnson and	T. Zhang, \textquotedblleft{}Accelerating Stochastic Gradient Descent Using Predictive Variance Reduction,\textquotedblright in \textit{Proc. Adv. Neural Inf. Process. Syst.},  vol. 1, pp. 315 - 323, 2013.
	
	\bibitem{Liu1989}
	D. C. Liu and J. Nocedal, \textquotedblleft{}On the limited memory BFGS method for large scale optimization,\textquotedblright \textit{Math. Programming}, vol. 45, no. 1,  pp. 503 - 528, Aug., 1989.
	
	\bibitem{pmlr-v51-moritz16}
	P. Moritz and R. Nishihara and M. I. Jordan, \textquotedblleft{}A Linearly-Convergent Stochastic L-BFGS Algorithm,\textquotedblright in \textit{Proc. 19th Int. Conf. Artif. Intell. Statist.}, vol. 51,   pp. 249 - 258, 2016.
	
	\bibitem{VIWMM}
	J. Taghia and A. Leijon, \textquotedblleft{}Variational Inference for Watson Mixture Model,\textquotedblright \textit{IEEE Trans. Pattern Anal. Mach. Intell.}, vol. 38, no. 9,   pp. 1886 - 1900, Sept., 2016.
	
	\bibitem{ARCG}
	A. Honkela,	T. Raiko,	M. Kuusela,	M. Torni	and Juha Karhunen , \textquotedblleft{}Approximate Riemannian Conjugate Gradient Learning for Fixed-Form Variational Bayes,\textquotedblright \textit{J. Mach. Learn. Res.},   vol. 11,  pp. 3235 - 3268, Dec., 2010.
	
	\bibitem{AmariShun-Ichi2016IGaI}
	S. Amari, \textit{Information Geometry and Its Applications},   Springer Japan, 2016.
	
	\bibitem{CVX}
	S. Boyd, and L. Vandenberghe,  \textit{Convex Optimization},   Cambridge University Press New York., 2004.
	
	\bibitem{scene}
	M.R. Boutell, J. Luo, X. Shen, and C.M. Brown, \textquotedblleft{}Learning multi-label scene classiffication,\textquotedblright \textit{Pattern Recognition},    vol. 37, no. 9, pp. 1757-1771, 2004.
	
	\bibitem{saa}
	B. T. Polyak and A. B. Juditsky, \textquotedblleft{}Acceleration of stochastic approximation by averaging,\textquotedblright \textit{SIAM J. Control Optim.},    vol. 30, no. 4, pp. 838-855, 1992.
	
	\bibitem{adam}
	D. Kingma and J. Ba,  \textquotedblleft{Adam: A method for stochastic optimization}\textquotedblright,   in \textit{3rd International Conference for Learning Representations, 2015}.
	
	
	\bibitem{Elec}
	V. Arzamasov, K. B\"{o}hm and P. Jochem,  \textquotedblleft{Towards Concise Models of Grid Stability
		}\textquotedblright, \textit{IEEE International Conference on Communications, Control, and Computing Technologies for Smart Grids}, Oct. 9-31,  2018.
	
	
	\bibitem{wifi1}
	R. Bhatt,  \textquotedblleft{Fuzzy-Rough Approaches for Pattern Classification: Hybrid measures, Mathematical analysis, Feature selection algorithms, Decision tree algorithms, Neural learning, and Applications}\textquotedblright, Amazon Books.
	\bibitem{wifi2}
	J. Rohra, B. Perumal, S. Narayanan, P. Thakur and R. Bhatt,  \textquotedblleft{User Localization in an Indoor Environment Using Fuzzy Hybrid of Particle Swarm Optimization \& Gravitational Search Algorithm with Neural Networks}\textquotedblright,  \textit{in Proceedings of Sixth International Conference on Soft Computing for Problem Solving}, pp. 286-295, 2017.
	
	
	
	\bibitem{IONO}
	V. G. Sigillito, S. P. Wing,  L. V. Hutton,  and  K. B. Baker,  \textquotedblleft{Classification of radar returns from the ionosphere using neural networks. Johns Hopkins APL Technical Digest}\textquotedblright, Johns Hopkins APL Technical Digest, 10, 262-266, 1989.
	
	
	\bibitem{BNA}
	V. Lohweg and H. Doerksen,  \textquotedblleft{Banknote authentication data set}\textquotedblright,   submitted.
	
	\bibitem{AsynADMM}
	R. Zhang and T. Kwok, \textquotedblleft{Asynchronous distributed ADMM for consensus optimization}\textquotedblright, in \textit{Proc. of the 31st  Int. Conf. Mach. Learn}, vol. 32, pp. 1701-1709, June, 2014.
	
	\bibitem{AsynBlockADMM}
	R. Zhu, D. Niu and Z. Li, \textquotedblleft{A Block-wise, Asynchronous and Distributed ADMM Algorithm for General Form Consensus Optimization}\textquotedblright, in \textit{arXiv:1802.08882}, Feb. 2018. 
	
	\bibitem{R48}
	T. G. Dietterich, \textquotedblleft{Approximate statistical tests for comparing supervised classification learning algorithms}\textquotedblright,  \textit{Neural Computation}, vol. 10, pp. 1895–1924, Oct. 1998.
	
	\bibitem{R49}
	J. Dem\v{s}ar, \textquotedblleft{Statistical Comparisons of Classifiers over Multiple Data Sets}\textquotedblright, \textit{J. Mach. Learn. Res.},  vol. 7, pp. 1- 30, Jan. 2006.
	
	
	\bibitem{MCHEN}
	M. Chen, B. Amos, L. Watson, J. Tyson, Y. Cao,
	C. Shaffer, M. Trosset, C. Oguz, and G. Kakoti, \textquotedblleft{}Quasi-Newton Stochastic Optimization Algorithm
	for Parameter Estimation of a Stochastic Model
	of the Budding Yeast Cell Cycle,\textquotedblright \textit{ IEEE/ACM Trans. Comput. Biol. Bioinf.}, vol. 16, no. 1, pp. 301-311, Nov. 2017.
	
	\bibitem{SHUANG}
	S. Huang, Y. Sun, and Q. Wu, \textquotedblleft{}Stochastic Economic Dispatch With Wind Using
	Versatile Probability Distribution and L-BFGS-B
	Based Dual Decomposition,\textquotedblright \textit{ IEEE Trans. Power Syst.}, vol. 33, no. 6, pp. 6254-6263, Nov. 2018.
	
	\bibitem{JRAFATI}
	J. Rafati and R. Marcia, \textquotedblleft{}Improving L-BFGS Initialization For
	Trust-Region Methods In Deep Learning,\textquotedblright in \textit{17th IEEE Int. Conf. Mach. Learn. App.}, Dec. 27-20, 2018.
	
	\bibitem{SCA}
	S. Scardapane and P. Lorenzo, \textquotedblleft{}Stochastic Training of Neural Networks via
	Successive Convex Approximations,\textquotedblright  \textit{IEEE Trans. Neural Netw. Learn. Syst.}, vol. 29, no. 10, Oct. 2018.
	
	
	\bibitem{AJA}
	A. Jalilzadeh, A. Nedi\'{c}, U. Shanbhag and F. Yousefian, \textquotedblleft{}A Variable Sample-size Stochastic Quasi-Newton Method
	for Smooth and Nonsmooth Stochastic Convex Optimization,\textquotedblright  \textit{ IEEE Conference on Decision and Control}, Dec. 17-19 2018.
	
	
\end{thebibliography}
\end{document}